\documentclass[10pt,a4paper]{article}
\newcommand{\rmodel}{~(\ref{eq:modello})}
\newcommand{\amodel}{~(\ref{eq:mode:abstract})}
\newcommand{\Dom}{{\rm Dom}}

\newtheorem{Theorem}{Theorem}

\newtheorem{Lemma}[Theorem]{Lemma}

\newtheorem{Remark}[Theorem]{Remark}
\newtheorem{Definition}[Theorem]{Definition}
\newtheorem{Example}[Theorem]{Example}
\newtheorem{Assumption}[Theorem]{Assumption}

\newcommand{\zaa}{\alpha}
 
\newcommand{\zg}{\gamma}
\newcommand{\ZDE}{\delta}
 \newcommand{\ZEP}{\epsilon}
\newcommand{\ZSI}{\sigma}
\newcommand{\zl}{\lambda} 
 \newcommand{\ZOM}{\omega}
 \newcommand{\ZOMq}{\Omega} 
 \newcommand{\zthe}{\theta}
\newcommand{\zt}{\tau}

\newcommand{\zzr}{\mathbb{R}}

 \newcommand{\zreal}{\Re{\textstyle e}\,} 
 
\newcommand{\intT}{\int_0^T}
\newcommand{\intt}{\int_0^t}

\newcommand{\zdiaform}{\mbox{~~\zdia}}

 \newcommand{\ZR}{\rangle}
\newcommand{\ZL}{\langle}
 \newcommand{\zdia}{~~\rule{1mm}{2mm}\par\medskip}
 
\newcommand{\ZIN}{\infty}
\newcommand{\zProof}{{\bf\underbar{Proof}.}\ }
 \newcommand{\ZD}{\;\mbox{\rm d}}
 
 \newcommand{\ZLA}{\label}

 \usepackage{epsfig}
\usepackage{graphics}
\usepackage[utf8]{inputenc}
\usepackage[english]{babel}
\usepackage{amsmath}
\usepackage{amsfonts}
\usepackage{amssymb}
\usepackage{makeidx}
\usepackage{graphicx}
\author{
L. Pandolfi\thanks{Retired from the Dipartimento di Scienze Matematiche ``Giuseppe Luigi Lagrange'', Politecnico di Torino, Corso Duca degli Abruzzi 24, 10129 Torino, Italy (luciano.pandolfi@formerfaculty.polito.it)}
}
\title{Controllability properties for   equations with memory of fractional type\thanks{
This papers fits into the research program of the GNAMPA-INDAM.}}
\begin{document}
 
 \maketitle 
 
 \begin{abstract} 
We study a general class of control systems with memory, which in particular includes systems with fractional derivatives and integrals and also the standard heat equation. We prove that the approximate controllability property of the heat equation is inherited by  every system   with memory in this class while controllability to zero is a singular property, which holds solely in the special case that the system   indeed reduces to the standard heat equation.
\end{abstract} 
 \section{Introduction}
In this paper we study the controllability properties of a general model for thermodynamical processes which in particular contains the model introduced by Coleman and Gurtin in~\cite{ColemanGurtinZEITSCH1967}:
\begin{equation}\ZLA{eq:modello}
\intt K(t-s)\frac{\ZD}{\ZD s} w(s)\ZD s=\Delta w(t)+\intt N(t-s)\Delta w(s)\ZD s +F(t)\,,\qquad w(0)=w_0\,.
 \end{equation}
 Here $\Delta$ is the laplacian and  $w=w(t)=w(x,t)$, $F=F(t)=F(x,t)$ with $x\in\ZOMq$ (a bounded region with $C^2$ boundary). The boundary conditions are described below in~(\ref{eq:BoundaryCONDI}). 
 
  We recall the notations   $*$ for the convolution, as in
 \[
K*w=\intt K(t-s)w(s)\ZD s \,.
 \]
With this notation, Eq.\rmodel\ takes the form
\[
K*w'=\Delta w+N*\Delta w+F\,,\qquad w(0)=w_0\,.
\]

 Systems of this type are widely studied (see the references below) and the focus in this paper is on the case that $K(t)$ and or $N(t)$ have a weak singularity at $0$.  The assumptions on $K$ and  $N$ used in this paper are described below. 
  The prototype are kernels of the form ($\Gamma$ is the Euler gamma function)
 \[
K(t)=\frac{1}{\Gamma(1-\zaa)} t^{-\zaa}\,,\qquad N(t)=\frac{1}{\Gamma(\zg)}t^{-(1-\zg)}\,.
 \]
 In this case Eq.\rmodel\ is the equation
\[
D^\zaa w=\Delta w+ J^\zg  w+F
\]
where $J^\zg$ is the Riemann-Liouville fractional integral and $D^\zaa$ the Caputo fractional derivative. As we don't need specific properties of fractional  integrals or derivatives we refer interested readers to~\cite{MAINARDIbook10FRACT}.

We do not exclude the case that either $K(t)$ is the Dirac  delta or the case $N=0$ i.e. the equations
\[
\begin{array}{l} {\rm either}\ 
w'=\Delta w+\intt N(t-s) w(s)\ZD s
\\
{\rm or}\  \intt K(t-s)w'(s)\ZD s= \Delta w+F\,.
\end{array}
\]
 
 The case $K(t)=\delta(t)$, $N(t)=0$ is the \emph{associated heat equation} to Eq.\rmodel: 
 \[u'=\Delta u+F\,,\qquad u(0)=w_0\]
 with the same boundary conditions in~(\ref{eq:BoundaryCONDI}) below as Eq.\rmodel.

The boundary conditions  are as follows. We fix a  nonempty relatively open subset $\Gamma_a$ of $\partial\ZOMq$. $\Gamma_a$ is the \emph{active part} of the boundary and we impose
\begin{equation}\ZLA{eq:BoundaryCONDI}
\mbox{  on $\Gamma_a$:}\  w=f={\rm boundary\ control;} \ \ \mbox{on $ \partial\ZOMq\setminus \Gamma_a$:} \   w=0 \,.
\end{equation}
So, when $w$
represents temperature, we are controlling the temperature on a part of the boundary.  

When the memory kernels $K(t)$ and $N(t)$ satisfy the assumptions~\ref{itemAssu}, system\rmodel\ is well posed, as specified in theorems~\ref{teo:EsiEVOLOPER} and~\ref{teo:EsiSOLUaffine}, i.e. for every initial conditions $w_0$ and every $F\in L^2\left (0,T;L^2(\ZOMq)\right )$ and $f\in L^2\left (0,T;L^2(\Gamma_a)\right )$ there exists a unique solutions (in the spaces specified in the theorems) which we denote $w_f(t;w_0,F)$.

\begin{Remark}[On the notations]
{\rm 
When one of the data $w_0$, $F$ or $f$ is zero,   it is not listed in the symbol used for the solution; i.e. $w(t;w_0)$ is the solution with $f=0$, $F=0$. When writing $w(t;w_0,F)$ we assume $f=0$ while   $w_f(t )$ denotes the solution with $w_0=0$, $F=0$.\zdia
}
\end{Remark}

When we study controllability under boundary controls we   assume $F=0$. Alternatively, we shall also study controllability under the action of the distributed control $F$  localized to an 
\emph{active subregion $\ZOMq_a$} of $\ZOMq$ (i.e. $F=0$ on $\ZOMq\setminus \ZOMq_a$). In this case we assume $f=0$.  \emph{The controllability properties of the associated heat equation are well known and  the goal of this paper is to understand at what extent these properties are inherited by the equation with memory~\rmodel.   }

Now we introduce the following definitions (we recall the notation $\mathcal{D}$ for the functions of class $C^\ZIN$ with compact support in the specified set):

 \begin{Definition}
 {\rm
 The reachable sets at time $T$ with controls localized respectively in $\ZOMq_a$ or $\Gamma_a$ are respectively  the following subsets  of $L^2(\ZOMq)$ 
 \begin{align*}
 &\mathcal{R}_{ d}(w_0;T)=\left \{ w(T;w_0,F)\,,\qquad F\in \mathcal{D}\left (\ZOMq_a\times(0,T)\right ) \right \}\\
&\mathcal{R}_{  b}(w_0;T)=\left \{ w_f(T;w_0 )\,,\qquad f\in \mathcal{D}\left (\Gamma_a\times(0,T)\right )\right \}\,.
 \end{align*}
 The system is approximately controllable at time $T$ (under distributed or boundary controls) when $\mathcal{R}_{ d}(0;T)$, respectively $\mathcal{R}_{  b}(0;T)$, is dense in $L^2(\ZOMq)$.
 
 The initial condition $w_0$ is controllable to hit the target zero at time $T$ when, respectively, $0\in \mathcal{R}_{ d}(w_0;T)$ or $0\in \mathcal{R}_{ b}(w_0;T)$.\zdia
 }
 \end{Definition}
 
As proved in~\cite{Lebeau95COMMpde}, \emph{the associated heat equation is approximately controllable at any time $T>0$ and every initial condition $w_0\in L^2(\ZOMq)$ can be steered to hit the target zero in any time $T>0$} (in both the cases, the controls are localized in any open subset of $\ZOMq_a$ or, respectively, in any relatively open subset $\Gamma_a$). The main result we are going to prove is:
\begin{Theorem}\ZLA{TeoINTRIUZcontroll}
We fix $\ZOMq_a$ and $\Gamma_a$. The set 
  $\ZOMq_a$ is an open subset of $\ZOMq$  while $\Gamma_a$ is a relatively open subset of $\Gamma=\partial\ZOMq$.
System\rmodel\ with the boundary condition~(\ref{eq:BoundaryCONDI}) has the following properties:
\begin{enumerate}
\item\ZLA{item:APPRX:TeoINTRIUZcontroll} For every $T>0$  the sets $\mathcal{R}_{ d}(w_0;T)$, respectively $\mathcal{R}_{ b}(w_0;T)$, are dense in $L^2(\ZOMq)$; 
\item\ZLA{item:NULL:TeoINTRIUZcontroll}   If $0\in \mathcal{R}_{ d}(w_0;T)$ or if $0\in \mathcal{R}_{ b}(w_0;T)$ for every $T>0$ (even for every $T$ in an open subset of $(0,+\ZIN)$) then system\rmodel\ is reducible to the associated heat equation.
\end{enumerate}
\end{Theorem}

The sense in which system\rmodel\  is reducible to the associated heat equation is specified in the lemmas~\ref{Lemma:NOaccumPOINT} and~\ref{Lemma:ReducASSOCheatEQ} and definition~\ref{defin:ReducHEAt}.

We sum up the main results of this paper: \emph{approximate controllability is inherited by system with memory while controllability to the target zero is a singular property, which holds solely for the associated heat equation.} These facts  have to be compared with the general results in~\cite{PANDOLFI1989LAA} where it is proved that both approximate controllability and controllability to the target zero (even controllability to zero in the state space sense) are not robust properties.

 The organization of the paper is as follows: the assumptions are presented in Sect.~\ref{Sect:Assump} (in the subsection~\ref{sebs:DiscuASSUM} it is shown that they are satisfied in particular by systems with fractional derivatives and integrals).
Approximate controllability is studied in Sect.~\ref{sect:controllabPROOFSapproximate} while controllability to the target zero is discussed in Section~\ref{sec:LacCONTRoZERO}. Ancillary results on the well posedness, stated in Section~\ref{Sect:Assump}, are proved in the Appendix~\ref{appe:ProofsEXI}.

\subsection{\ZLA{Sect:Assump}The assumptions and the definitions of the solutions}
 
We relay on Laplace transform techniques to study Eq.\rmodel. So, the assumptions are expressed in the frequency domain.    We introduce the following notation: if $\zaa\in (0,\pi)$, $\Sigma_{\zaa}$ is the sector of the complex plane
 \begin{equation}\ZLA{eq:GeneDEFIsector}
 \Sigma_{\zaa}=\{\zl\,:\ |{\rm Arg}\zl|<\zaa\,,\quad \zl\neq 0\}\,.
 \end{equation}

Let $A$ be the operator   in $L^2(\ZOMq)$:
\[
\Dom\,A=H^2(\ZOMq)\cap H^1_0(\ZOMq)\,,\qquad A w=\Delta w.
\]
   
The following properties are known:
\begin{enumerate}
\item the operator $A$ is selfadjoint positive with compact resolvent;
\item $\sup\ZSI(A)<0$;
\item Let  $\zthe_A\in (0,\pi/2)$. There exist $\ZOM>0$ and   $M>0$ (which depends on $\zthe_A$ and $\ZOM$)    such that the following inequalities hold for every $\zl\in \Sigma_{\zthe_A+\pi/2}$:
\begin{equation}\ZLA{eq:DiseGeneSGolo}
\left \|(\zl I-A)^{-1}\right \|\leq\left\{\begin{array}{l}\frac{M}{|\zl+\ZOM|}\\[2mm] 
\leq  \frac{M}{|\zl|} 
\end{array}\right.
\end{equation}
(here and below we use $M$ to denote a generic constant, not the same at every occurrence).

\emph{We intend that $\zthe_A$ and $\ZOM$ have been chosen and kept fixed in the paper.}

\item the operator $A$ generates a holomorphic semigroup $e^{At}$, which admits holomorphic extension to the sector $\Sigma_{\zthe_A}$.
\end{enumerate}

We introduce the Green operator $G$ which is defined as follows:
\[
u=Gf\qquad \iff \quad \left\{\begin{array}{l}
\Delta u=0\\
u=f \ {\rm on}\ \Gamma_a\,;\quad u=0\ {\rm on}\ \partial\ZOMq\setminus\Gamma_a\,.
\end{array}\right.
\]
It is known that there exists $\ZSI_0\in (0,1)$ (any $\ZSI_0<1/4$ but the precise value is not needed. Hence most of  the arguments in this paper can be applied also to different boundary conditions, when the exponent takes different values) such that ${\rm im}\, G\subseteq \Dom\, (-A)^{\ZSI_0}$ so that from the interpolation inequality of the fractional powers of $A$ we get:
\begin{equation}
\ZLA{eq:DiseqPerG}
\left \|(\zl I-A)^{-1}AG\right |\leq \dfrac{M}{|\zl+\ZOM|^{\ZSI_0}}\leq \frac{M}{|\zl|^{\ZSI_0}}\quad \forall\zl\in\Sigma_{\zthe_A+\pi/2}.
\end{equation}

It is now well understood that the operator $G$ can be used to insert the boundary control in the equation\rmodel\ (see~\cite{LasieckaTriggianiLIBROencVOL1,PPZ}). 
We use the following equality, which is correct in $\left  (\Dom \,A\right )'  $:
\[
\Delta w=\Delta (w-Gf)=A(w-Gf)=Aw-AGf
\]
and we rewrite
Eq.\rmodel\   in $\left  (\Dom \,A\right )'  $ as follows:
\begin{equation}\ZLA{eq:mode:abstract}
K*w'= Aw+N*Aw+AG\underbrace{\left ( f+N*f  \right )}_{g}+F\,,\qquad w(0)=w_0\,.
\end{equation}
The solutions of this equation have to be found    in $\left  (\Dom \,A\right )'  $ but, as we shall see, when $w_0\in L^2(\ZOMq)$ we can give a sense to the solutions in the space $L^2(\ZOMq)$.

In order to understand the assumptions on the memory kernels, let us compute formally the Laplace transform   of\amodel. We get
\begin{multline*}
\hat K(\zl)\left (\zl \hat w(\zl)-w_0\right )
=J(\zl) A\hat w(\zl)-J(\zl)AG\hat f(\zl)+\hat F(\zl)\,,\qquad J(\zl)=1+\hat N(\zl)
\end{multline*}
and so:
\begin{multline}
\ZLA{eq:TrasfoRISOLTA}
\hat w(\zl)=\left (\zl\hat K(\zl)I-J(\zl)A\right )^{-1}\left [ \hat K(\zl) w_0-J(\zl)AG\hat f(\zl)+\hat F(\zl)\right ]\\
=\left (\zl\hat K(\zl)I-J(\zl)A\right )^{-1}\hat K(\zl)w_0+\left (\zl\hat K(\zl)I-J(\zl)A\right )^{-1}\hat F(\zl)\\
-
\left (\frac{\zl\hat K(\zl)}{J(\zl)}-A\right )^{-1}AG\hat f(\zl)\\
=\left (\zl\hat K(\zl)I-J(\zl)A\right )^{-1}\hat K(\zl)w_0+\left (\zl\hat K(\zl)I-J(\zl)A\right )^{-1}\hat F(\zl)\\
+
\left [G\hat f(\zl)-\zl\hat K(\zl) \left (\zl\hat K(\zl)I-J(\zl)A\right )^{-1}G\hat f(\zl)  \right ]\,.
\end{multline}
Suggested by this expression we assume:
\begin{Assumption}[On the memory kernels]
{\rm
{~}

\ZLA{itemAssu}
\begin{enumerate}
\item\ZLA{itemAssuKERNEL00} 
$K(t)=k_0\ZDE(t)+K_0(t)$ where $k_0\geq 0$ and $\ZDE(t)$ is the Dirac  delta and    both the kernels $K_0(t)$ and $N(t)$ are real valued and integrable on $(0,+\ZIN)$.  Hence, the Laplace transforms $\hat K(\zl)$ and $\hat N(\zl)$ have the following properties
\begin{itemize}
\item the Laplace transforms of  $\hat K(\zl)$ and $\hat N(\zl)$ exist and are holomorphic functions on $\zreal\zl>0$.
\item 
$\lim _{\stackrel{|\zl|\to+\ZIN}{\zreal\zl>0 }} \hat K(\zl)=k_0\geq 0  $.
\item $\lim _{\stackrel{|\zl|\to+\ZIN}{\zreal\zl>0 }} \hat N(\zl)= 0$.
\end{itemize}

\item\ZLA{itemAssuKERNEL2}  There exists $\zthe\in (0,\pi/2)$ such that $\hat K(\zl)$ and $\hat N(\zl)$ (originally defined for $\zreal\zl>0$) admit holomorphic extensions to $\Sigma_{\zthe+\pi/2}$.

 \emph{We recall that $\zthe_A\in (0,\pi/2)$ has been chosen once and for all,  and kept fixed in the paper. Hence,  it is not restrictive to assume $0<\zthe<\zthe_A$.}
 \item\ZLA{itemAssuKERNEL3}  $J(\zl)=1+\hat N(\zl)\neq 0$ for $\zl\in\Sigma_{\zthe +\pi/2}$.
\item\ZLA{itemAssuKERNEL4}  the map
 \[\zl
 \mapsto \zl\hat K(\zl)/J(\zl) 
 \]
 transforms $\Sigma_{\zthe+\pi/2} $ to $\Sigma_{\zthe_A+\pi/2}$.
\item\ZLA{itemAssuKERNEL2ADDE}
Assumptions on the asymptotics of   $\hat K(\zl)$ and $\hat N(\zl)$:
\begin{enumerate}

\item \ZLA{itemAssuASYMPkernelK} asymptotics of   $\hat K(\zl)$:
\begin{enumerate}
 \item \ZLA{itemAssuKERNEL2ADDE0}
$\lim _{\stackrel{|\zl|\to+\ZIN}{\zl\in \Sigma_{\zthe+\pi/2}}} \hat K(\zl)=k_0\geq 0  $\,.
\item \ZLA{itemAssuKERNEL2ADDE1}
 there exists   positive numbers $M$ and    $R$ and an exponent  $\zg_0\in (0,1]$ such that $|\zl \hat K(\zl)|>M|\zl|^{\zg_0}$ for $\zl\in \Sigma _{\theta+\pi/2}$, $|\zl|>R$;
\item \ZLA{itemAssuKERNEL2ADDE2}
  $\lim _{\stackrel{\zl\to 0}{\zl\in \Sigma_{\zthe+\pi/2}}}\zl \hat K(\zl)=0$.
 
 \end{enumerate}
\item\ZLA{itemAssuKERNEL2ADDE0nuclN}  asymptotics of $\hat N(\zl)$:
 $ \lim _{\stackrel{|\zl|\to+\ZIN}{\zl\in \Sigma_{\zthe+\pi/2}}} \hat N(\zl)=0 $ so that there exists $r>0$ such that when $\zl\in \Sigma_{\zthe+\pi/2}$ and $|\zl|>r$ we have
\[
\frac{1}{2}< |J(\zl)|=|1+\hat N(\zl)|<\frac{3}{2}. \zdiaform
\]

\end{enumerate}

\end{enumerate}
}
\end{Assumption}

The fact that these assumptions are satisfied in significant cases is seen in Sect.~\ref{sebs:DiscuASSUM}.
In particular, note that the constant  $k_0 $ in item~\ref{itemAssuKERNEL2ADDE0}   is strictly positive if $K(t)$ is 
the Dirac  delta, as in the case of the associated heat equation.

\begin{Remark}\ZLA{rema:sulNONnullo}
{\rm 
We note:
\begin{itemize}
\item We are interested on controllability on bounded intervals $[0,T]$. So, the assumption that $K(t)$ and $N(t)$ are integrable on $(0,+\ZIN)$ is in fact the conditions that they are integrable on $(0,T)$ and then extended with zero on $(T,+\ZIN)$;
\item $\hat K(\zl)\neq  0$ in $\Sigma_{\zthe+\pi/2}$ since $\zl\hat K(\zl)/J(\zl)\in \Sigma_{\zthe_A+\pi/2}$ and $0\notin \Sigma_{\zthe_A+\pi/2}$ (see the definition~(\ref{eq:GeneDEFIsector}));
\item Let $\ZOM$ be the constant  in~(\ref{eq:DiseGeneSGolo}), hence $\ZOM>0$. Then $\zl\hat K(\zl)+\ZOM J(\zl)\neq 0$ in $\Sigma_{\zthe+\pi/2}$ since otherwise there should exists $\zl\in  \Sigma_{\zthe+\pi/2}$  such that $\zl\hat K(\zl)/J(\zl)=-\ZOM\notin \Sigma_{\zthe_A+\pi/2}$.
\item the assumptions in item~\ref{itemAssuKERNEL2ADDE} imply 
\[
\lim _{|\zl|\to+\ZIN,\ \zl\in \Sigma_{\zthe+\pi/2}} \frac{\zl \hat K(\zl)}{J(\zl) }  =k_0\,.\zdiaform
\]
\end{itemize}
}
\end{Remark}
 
Using~(\ref{eq:DiseGeneSGolo}) and~(\ref{eq:DiseqPerG}) we see that the previous assumptions imply the following inequalities  (with possibly different values of the constants denoted $M$). These inequalities hold on the sector  $\Sigma_{\zthe+\pi/2}$:
\begin{equation}\ZLA{eq:ineqMODIFICATE}
\left\{\begin{array}{l}\displaystyle
\left \| \left (  \zl\hat K(\zl) I-J(\zl)A\right )^{-1}\right \|\leq
\left\{\begin{array}{l}
 \frac{M}{\left 	|\zl\hat K(\zl)+\ZOM J(\zl)\right |}
 \\[2mm] 
 M \frac{1}{|\zl\hat K(\zl)|}
 \end{array}\right.\\
 \\
 \displaystyle
 \left \| \left (\frac{\zl\hat K(\zl)}{J(\zl)}I-A\right )^{-1}AG\right\|\leq 
 \left\{\begin{array}{l}
 M\left | \frac{\zl\hat K(\zl)}{J(\zl)}+\ZOM\right |^{-\ZSI_0}\\[2mm]
 M\frac{|J(\zl)|^{\ZSI_0}}{|\zl \hat K(\zl)|^{\ZSI_0}} \,.
 \end{array}\right.
 \end{array}\right.
\end{equation}
  
The  main results we prove in this paper   depend  on the definition and properties of the solutions of Eq.\rmodel, i.e.\amodel,  collected in the following theorems:
\begin{Theorem}\ZLA{teo:EsiEVOLOPER}
There exists an \emph{evolution operator $E(t)$} i.e. a function from $t\in[0,+\ZIN)$ to $\mathcal{L}\left (L^2(\ZOMq)\right )$ with the following properties:
\begin{enumerate}
\item\ZLA{teo:EsiEVOLOPERitem1} $E(t)$ is   continuous for $t>0$,  bounded on $[0,+\ZIN)$ and $\lim _{t\to 0}E(t)w_0=w_0$ for every $w_0\in L^2(\ZOMq)$;
\item\ZLA{teo:EsiEVOLOPERitem3} the functions $E(t)$ admits a holomorphic extension to the  sector $\Sigma_{\zthe }$ (the extension is denoted $E(z)$ and $\zthe $ is the angle in item~\ref{itemAssuKERNEL2} of the set of the assumptions~\ref{itemAssu}) and for each $z\in\Sigma_\zthe$  (hence not for $z=0$) we have $E(z) w\in\Dom\, A$ for every $w\in L^2(\ZOMq)$ (we recall the nonrestrictive assumpotioin $0<\zthe<\zthe_A$);
\item\ZLA{teo:EsiEVOLOPERitem4} the function $AE(z)$ is holomorphic on $\Sigma_\zthe$;
\item  \ZLA{teo:EsiEVOLOPERitem2} the Laplace transform of $E(t)w_0 $ is $\left (\zl\hat K(\zl)I-J(\zl)A\right )^{-1} \hat K(\zl) w_0$.
\end{enumerate}
 
\end{Theorem}

The last property in theorem~\ref{teo:EsiEVOLOPER} justifies the following definition:
\begin{Definition}
{\rm
A mild solution of system\amodel\ (i.e.\rmodel) when $F=0$, $f=0$ is $w(t;w_0)=E(t)w_0$ for every initial condition $w_0\in L^2(\ZOMq)$. 
}
\end{Definition}

We recall the notation   $w_f(t;w_0,F)$ for the solution, and we recall that when $w_0$, $ F$ or $f$  is zero, the corresponding symbol is omitted.

The following theorem justify the definition~\ref{DefiSOLU} below of the solutions $w(t;F)$ and $w_f(t)$:

\begin{Theorem}\ZLA{teo:EsiSOLUaffine}
The following properties hold:
\begin{enumerate}
\item\ZLA{teo:EsiSOLUaffineITEM1} For every $T>0$ there exists a linear continuous operator  operator $\mathcal{E}_{d,T}$ from $ L^2\left (0,T;L^2(\ZOMq)\right )$ to itself such that
\begin{enumerate}
\item if $T_1>T$ and $F\in L^2\left (0,T_1;L^2(\ZOMq)\right )$ then $\left (\mathcal{E}_{d,T_1} F\right )_{|_{(0,T)}} =\mathcal{E}_{d,T} \left ( F _{|_{(0,T)}} \right ) $. So, we can define  $\left (\mathcal{E}_{d }F\right )(t)=\left (\mathcal{E}_{d,T }F\right )(t)$ with any $T>t$;
\item the Laplace transform of $\mathcal{E}_{d }F$ is defined for $\zreal\zl >0$ and it is equal to $\left (\zl\hat K(\zl)I-J(\zl)A\right )^{-1}\hat F(\zl)$.
\item\ZLA{teo:EsiSOLUaffine:itemCONdistri} if $F\in \mathcal{D}\left (\ZOMq\times(0,T)\right )$ then the function $t\mapsto  \left (\mathcal{E}_{d }F\right )(t)$ is continuous from $[0,T]$ to $ \Dom\, A^k $ for every $T>0$ and every $k\geq 0$.
\end{enumerate}
\item\ZLA{teo:EsiSOLUaffineITEM2} For every $T>0$ there exists   a linear continuous operator $\mathcal{E}_{b,T}$ from $L^2\left (0,T;L^2(\Gamma_a)\right ) $ to  $L^2\left (0,T;L^2(\ZOMq)\right ) $  such that
\begin{enumerate}
\item if $T_1>T$ and $f\in L^2\left (0,T_1;L^2(\Gamma_a)\right )$ then $\left (\mathcal{E}_{b,T_1} f\right )_{|_{(0,T)}} =\mathcal{E}_{b,T} \left ( f _{|_{(0,T)}} \right ) $. So, we can define  $\left (\mathcal{E}_{b }f\right )(t)=\left (\mathcal{E}_{b,T }f\right )(t)$  with any $T>t$;
\item the Laplace transform of $\mathcal{E}_{b }f$ is defined for $\zreal\zl >0$ and it is equal to $-\left (\zl\hat K(\zl)I-J(\zl)A\right )^{-1}J(\zl)AG \hat f(\zl)$.
\item\ZLA{teo:EsiSOLUaffine:itemCONbound} if $f\in \mathcal{D}\left (\Gamma_a\times(0,T)\right )$ then the function $t\mapsto  \left (\mathcal{E}_{b }f\right )(t)$ is continuous from $[0,T]$ to   $ L^2(\ZOMq) $ for every $T>0$.
\end{enumerate}
\end{enumerate}
\end{Theorem}

These results justify:
\begin{Definition}\ZLA{DefiSOLU}
{\rm 
When $w_0=0$, $f=0$ we define $w(\cdot, F)=\mathcal{E}_{d }F$; when $w_0=0$, $F=0$ we define
$w_f(\cdot )=\mathcal{E}_{b }f$ and, in general, $w_f(\cdot;w_0,F)=E(t) w_0+\mathcal{E}_{d }F+\mathcal{E}_{b }f$. 
}
\end{Definition}

The proofs of the previous theorems~\ref{teo:EsiEVOLOPER} and~\ref{teo:EsiSOLUaffine} follow  the lines already used in several special cases and for completeness they are sketched in the appendix~\ref{appe:ProofsEXI}.

\begin{Remark}
{\rm
It is easy to understand the importance that  the statement in item~\ref{teo:EsiSOLUaffine:itemCONbound} of theorem~\ref{teo:EsiSOLUaffine} has in the definition of the reachable set under boundary controls, since it is known that in general $w_f(t)$ is not continuous
if $f$ is solely square integrable,  not even in the case of 
the associated heat equation, as seen in the example~\cite[p.~217]{LIONSlibro68}. 
 Example~\ref{EXE:nOcONTIN} in the appendix~\ref{appe:ProofsEXI} shows that even the function $t\mapsto w (t;F)$ is not continuous when $K(t)$ is not the Dirac  delta, and this explains the reason why we stated explicitly   the result in item~\ref{teo:EsiSOLUaffine:itemCONdistri}.\zdia
}
\end{Remark}

We shall use also representation formulas for the operators $E(t)$, $\mathcal{E}_d(t)$ and $\mathcal{E}_b(t)$.  In the case of the associated heat equation, $\mathcal{E}_d(t)$ and $\mathcal{E}_b(t)$ are just convolutions with    $E(t)$ which in this case is $e^{At}$. This is not the case in general, since the decaying factor $K(\zl)$ in $\hat w(\zl;w_0)$ is not present in the expression of $\hat w(\zl;F)$ and $\hat w_f(\zl )$.

The representations, derived in the appendix~\ref{appe:ProofsEXI}, are as follows. Let 
the  path of integration $G_\ZEP$ be composed by the following two half lines and   circular arch (this is the usual integration path of the theory of holomorphic semigroups, see~\cite[Sect.~1.7]{Pazy}. 
The path is represented in Fig.~\ref{fig:path} of the appendix~\ref{appe:ProofsEXI}):
 \begin{equation}\ZLA{eq:intepath}
G_{\pm}\,:\ \zl=s\left [\cos\zaa\pm\sin\zaa\right ]\,\ s\in [\ZEP,+\ZIN)\,,\qquad  \zl=\ZEP e^{-\zt}\,,\ -\zaa<\zt<\zaa\,.
 \end{equation}
The angle $\zaa\in (\pi/2,\zthe)$ is fixed. Then we have
\begin{equation}\ZLA{eq:rAPPreDIwinidata}
w(t;w_0)=E(t)w_0=\frac{1}{	2\pi i}\int _{G_\ZEP} e^{t\zl} \hat K(\zl)\left (\zl\hat K(\zl)I-J(\zl)A\right )^{-1}w_0\ZD \zl\,.
\end{equation}
Of course the improper integral is computed as the limit  for $R\to+\ZIN$ of integrals on $G_{R,\ZEP}=G_{\ZEP}\cap \{\zl\,, \ |\zl|<R\}$.

The operators  $\mathcal{E}_d(t)$ and $\mathcal{E}_b(t)$ are given by the following convolutions:
\begin{equation}
\ZLA{rappreDISTRIBOUNDcontro}
w(t;F)=\intt  \mathcal {E}(t-s)F(s)\ZD s\,,\qquad w_f(t)=-\intt \mathcal {E}(t-s) AG g(s)\ZD s
\end{equation}
where
\begin{equation}
\ZLA{RappreEcorsivo}
\mathcal{E}(t)=\frac{1}{	2\pi i}\int _{G_\ZEP} e^{t\zl}  \left (\zl\hat K(\zl)I-J(\zl)A\right )^{-1} \ZD \zl 
\end{equation}
and $g$ is the function in~(\ref{eq:mode:abstract}).
As proved in the appendix~\ref{appe:ProofsEXI}, the integral converges uniformly for $t\in [a,T]$, $a>0$ and in $L^1\left (0,T; \mathcal{L}(X)\right )$ for every $T>0$.

\subsubsection{\ZLA{sebs:DiscuASSUM}Discussion of the assumptions}

The powers of $\zl$ denotes the principal value, so that $\arg\zl^\zg=\zg\arg\zl$ with $-\pi\leq \arg\zl<\pi$.

We recall $A=\Delta $ with domain $H^2(\ZOMq)\cap H^1_0(\ZOMq)\subseteq L^2(\ZOMq)$ so that $\theta_A$ is any number in $(0,\pi/2)$ and 
we discuss the assumptions on the memory kernels in the following significant cases:
\begin{enumerate}
\item\ZLA{SIGN1} $K(t)=\ZDE(t)$, $N(t)=\frac{1}{\Gamma(\zg)} t^{ \zg-1 }$ so that $\hat K(\zl)=1$, $\hat N(\zl)=\frac{1}{	\zl^\zg}$ where $\zg\in (0,1)$. Then
\[
\frac{\zl\hat K(\zl)}{J(\zl)	}=\frac{\zl^{1+\zg}}{1+\zl^{\zg}}\,;
\]
\item\ZLA{SIGN2}  $K(t)=\frac{1}{\Gamma(1-\zaa)} t^{-\zaa}$, $N(t)=0$ with $\zaa\in (0,1)$. Then
  $\zl \hat K(\zl)=\zl^{\zaa}$ and
\[
\frac{\zl\hat K(\zl)}{J(\zl)	}=\zl^{\zaa}\,;
\]
\item\ZLA{SIGN3} $K(t)=\frac{1}{\Gamma(1-\zaa)} t^{-\zaa}$, $N(t)=\frac{1}{\Gamma(\zg)} t^{ \zg-1}$
so that $\frac{\zl\hat K(\zl)}{ (1+\hat N(\zl)  )}=\zl^{1+\zaa}/\left (1+1/\zl^\zg\right )$ where $\zaa$ and $\zg$ belong to $(0,1)$. Then,
\[
\frac{\zl\hat K(\zl)}{J(\zl)	}=\frac{\zl^{ \zaa+\zg }}{1+\zl^{\zg}}\,;
\]
\end{enumerate}

  The conditions in the items~\ref{itemAssuKERNEL00}-\ref{itemAssuKERNEL3} and the asymptotic conditions of item~\ref{itemAssuKERNEL2ADDE} clearly hold in the three examples,  and we must prove the existence of a sector $\Sigma_{\zthe+\pi/2}$ over which the condition in item~\ref{itemAssuKERNEL4} hold, i.e. $\zl\hat K(\zl)/J(\zl)$ transforms $\Sigma_{\zthe+\pi/2}$ into a sector
$\Sigma_{\pi-\ZEP}$ for some $\ZEP>0$.
This assumption is clearly satisfied in   case~\ref{SIGN2} since $\zaa\in (0,1)$. 
We prove that it is satisfied in   case~\ref{SIGN1}. We consider   $\arg\zl\geq 0$ (similar computations when $\arg\zl\leq 0$   lead  to the same condition on $\ZEP$).  If $\arg\zl>0$ we have 
\[
\arg\frac{\zl^{1+\zg}}{1+\zl^\zg}=(1+\zg)\arg\zl- \arg (1+\zl^\zg  )<(1+\zg)\arg\zl\,.
\]
The required condition holds of we choose $\ZEP\in \left (0,(1-\zg)\pi/2\right )$.

Case~\ref{SIGN3} is more interesting since in this case the condition in item~\ref{itemAssuKERNEL4} imposes a link among the exponents $\zaa$ and $\zg$.  It is convenient to rename $\ZEP$ as $\pi\ZEP$ and $\zthe$ as $\zthe\pi/2$ so that now $\zthe\in (0,1)$.
Then, the  condition  in item~\ref{itemAssuKERNEL4} of the assumptions~\ref{itemAssu} 
can   be written as follows: there exists $\zthe\in (0,1)$ and $\ZEP>0$ such that both the following conditions hold when $-(\pi/2)(\zthe+1)\leq \arg\zl\leq (\pi/2)(\zthe+1)$:

\begin{equation}\ZLA{eq:SubsSIgnCAsesCase3}
 -\pi(1-\ZEP)<\arg\frac{\zl^{\zaa+\zg}}{1+\zl^\zg}<\pi(1-\ZEP)\quad 
\end{equation}
 
We examine the case $0\leq \arg\zl\leq \frac{\pi}{2}\left (\zthe+1\right )$. Computing the limit for $\zl\to 0$ along the line $\arg\zl=\left (\zthe+1\right )(\pi/2)$ we get  
\[
(\zaa+\zg)(\pi/2)(\zthe+1)<\pi (1-\ZEP)\quad \mbox{which implies the necessary condition $\zaa+\zg<2$ }\,.
\]

When $\zreal\zl\geq 0$ the left inequality in~(\ref{eq:SubsSIgnCAsesCase3}) is always satisfied since $\arg(1+\zl^\zg)\leq \zg\arg\zl$.
We prove that the inequality from above is satisfied, for a suitable $\zthe\in (0,1)$ and $\ZEP>0$, when $\zaa+\zg<2$. In fact,
\[
(\zaa+\zg)\arg\zl-\arg(1+\zl^\zg)\leq (\zaa+\zg)\arg\zl\leq (\zaa+\zg)\frac{\pi}{2}(\zthe+1) 
\]
and we want the existence of $\zthe_0>0$ such that the following holds for $\zthe<\zthe_0$:
\[
\frac{1}{2}(\zaa+\zg)(\zthe+1)\pi <\pi\,.
\]
This is achieved  provided that
\[
\zthe<\zthe_0\quad {\rm where}\quad \zthe_0=\frac{2-(\zaa+\zg)}{\zaa+\zg}>0\,.
\]

A similar analysis shows that when $\zaa+\zg<2$ the required condition~(\ref{eq:SubsSIgnCAsesCase3}) can be satisfied also in a sector $-(\pi/2)(\zthe+1)\leq\arg\zl\leq 0$.

\subsection{Comments on previous results}
The first instance  of system\rmodel\ that has been studied is the 
  Colemann-Gurtin equation, i.e. Eq.\rmodel\ 
when $K(t) $ is the Dirac  delta and of course the first results concern the case that the kernel $N(t)$ is smooth. In this case, the first study of controllability seems to be 
  the paper~\cite{BaumeisterJMAA1983} 
which proves \emph{lack of exact controllability.} This result is expected of course, since exact controllability does not hold for the associated heat equation. The proof of approximate controllability   in~\cite{BarbuDIFFINTEQ2000} is for kernels $K(t)=\ZDE(t)$ and 
  \[ 
  N(t)=\sum _{j=1}^n  a_j e^{-\zaa_j t} +\sum _{j=1}^m \int_{I_j} b_j(s) e^{-st }\ZD s
  \]
  where $ \zaa_j>0 $, $a_j\geq 0$;  $ b_k(s)>c_k>0 $ are integrable on the bounded intervals $ I_j $. Approximate controllability is then proved for  any $ N\in  H^1(0,T)$   (and $K(t)=\ZDE(t)$) in~\cite{HalanayJMAA2015,HalanayDCDS-A2015}. In this paper   lack of controllability to hit the target zero has also been proved for every $N\in H^1(0,T)$,     after the preliminary results in~\cite{GuerreroESAIM2013,HalanaySCL2012,HalanayJMAA2014}.
  
  Approximate controllability for a heat equation perturbed by a memory term $\intt H(t-s)u(s)\ZD s$ when the relaxation kernel has compact support is proved via Carleman inequalities   in~\cite{LavanyaBALACH}. A positive result on controllability to the target zero for this class of systems when $H(t)$ has  very special properties is proved in~\cite{TaoGAOcontroZEROjmma16} via Carleman estimates.
  
  Controllability of systems with fractional derivatives (in time) have been advocated in particular by M. Yamamoto, and we cite the paper~\cite{FujishiroYAMAMOTO2014} where approximate controllability is proved for the system~\rmodel\ with boundary controls, in the case that $N(t)$ is   the Dirac  delta and $K(t)=(1/\Gamma(1-\zaa)) t^{-\zaa} $ (see also~\cite{WarmaDISCRCONTDYNSYS2016,WARMAapplANAL2017}).

\section{\ZLA{sect:controllabPROOFSapproximate}Approximate controllability}

In this section we prove   the statement in item~\ref{item:APPRX:TeoINTRIUZcontroll} of theorem~\ref{TeoINTRIUZcontroll}, i.e. we prove that the system with memory\rmodel\ inherits the approximate controllability properties of the associated heat equation. In spite of the fact that distributed and boundary controls have special features, it is convenient to introduce suitable notations using which we can perform most of the computations in a unified way.

We recall that the distributed control $F$ is zero on $\ZOMq\setminus\ZOMq_a$. 
We introduce, for $F\in L^2\left (0,T;L^2(\ZOMq_a)\right )$,

\[
\left (B_0F\right )(x)=\left\{\begin{array}
{lll}
F(x)&{\rm if}& x\in\ZOMq_a\\
0&{\rm if}&x\in\ZOMq\setminus\ZOMq_a  
\end{array}\right.
\]
 (in this section it might be $\ZOMq_a=\ZOMq$. In this case, $B=I$). 
The control process\amodel\ is
\[
K*w'= Aw+N*Aw-AG\left (f+N*f\right )+B_0 F\,,\qquad w(0)=w_0\,.
\]
We can subsume both  the controls $f$ and $F$ in the same expression by writing
\begin{equation}\ZLA{eq:mode:unified}
K*w'= Aw+N*Aw+AB\left (g+N*g\right ) \,,\qquad w(0)=w_0 
\end{equation}
where ($R(t)$ is the resolvent kernel of $N(t)$, see below)
\begin{equation}\ZLA{eq:PerLAContraPPDefig}
\left\{\begin{array}{ll}{\rm either}\ B=G\ {\rm and}\ g=- f  \\
{\rm or} \ B=A^{-1}B_0\ {\rm and}\ g=F-R*F
\end{array}\right.
\end{equation}
and then we can study approximate controllability  of~(\ref{eq:mode:unified}) under the action of the control $g\in L^2(0,T;U)$ where now $U$ is either $L^2(\Gamma_a)$ or $L^2(\ZOMq_a)$.
The solution of~(\ref{eq:mode:unified}) with $w_0=0$ is denoted $w_g(\cdot)$.

Once this transformation has been done, the reachable set is denoted simply $\mathcal{R}(T)$:
\[
\mathcal{R}(T)=\left\{\begin{array}{l}
{\rm either}\ \mathcal{R}_d(0,T)\\
{\rm or}\ \mathcal{R}_b(0,T)\,.
\end{array}\right.
\]
The set $\mathcal{R}(T)$ is a subspace of $L^2(\ZOMq)$ and $\mathcal{R}(T_1)\subseteq \mathcal{R}(T_2)$ if $T_1<T_2$. Then we define
\[
\mathcal{R}(\ZIN)=\bigcup _{T>0}\mathcal{R}(T)\,.
\]

 Note that the ``real control'' in~(\ref{eq:mode:unified}) is not $g$ but $g+N*g$. Then it is convenient to keep in mind the following facts.

The resolvent kernel $R(t)$   of $N(t)$ is the unique solution of
\[
R+N*R=N\,.
\]
We recall that the transformation
\[
F\mapsto F+N*F
\]
is continuous and boundedly invertible in any space $L^2(0,T;U)$ (also $T=+\ZIN$ when $N(t)\in L^1(0,+\ZIN)$) and the inverse is 
\[
F\mapsto F-R*F\,.
\] 
 \emph{So, in the particular case of the distributed control we have $g+N*g=F$.} In paricular, if $F\in\mathcal{D}\left (\ZOMq_a\times(0,\ZIN)\right )$ then $g+N*g=F$ has this same regularity. Things are slightly different in the case of the boundary controls. In this case
 $f\in\mathcal{D}\left ( \Gamma_a \times (0,+\ZIN)\right )$ so that $N*f\in C^\ZIN\left ( \Gamma_a\times (0,+\ZIN)\right )$ and $N*f=0$ on a suitable interval $(0,\ZEP)$  ($\ZEP>0$ depends on $f$). 
 Furthermore:
 
 \begin{itemize}
 \item The set of the functions $g+N*g$ is dense in $L^2(0,T;U)$ both in the case $U=L^2(\ZOMq_a)$ and $U=L^2(\Gamma_a)$;
 \item
In both the cases, the Laplace transform $\hat g(\zl)$ decays faster then $1/|\zl|^k$ when $|\zl|\to+\ZIN$ in $\zreal\zl\geq 0$, for every $k>0$;
\item
The Laplace transform  of $R(t)$ is $\hat R(\zl)=\hat N(\zl)/J(\zl)$.
\end{itemize}

Now we recall that our goal is the proof that 
the subspace $\mathcal{R}(T)$ is dense in $L^2(\ZOMq)$ for every $T>0$.
The proof is in three steps:
\begin{enumerate}\ZLA{itemizeTHEthreeSTEPS}
\item\ZLA{itemizeTHEthreeSTEPSstep1}   we characterize the elements of $\mathcal{R}(T)^\perp$ for every $T>0$,    $T=+\ZIN$ included;
 \item\ZLA{itemizeTHEthreeSTEPSstep2}  we prove that $\mathcal{R}(\ZIN)$ is dense in $L^2(\ZOMq)$;
 \item\ZLA{itemizeTHEthreeSTEPSstep3} then we prove that $\mathcal{R}(T)$ is dense in $L^2(\ZOMq)$ at any time $T>0$.
\end{enumerate}

Note that in the proofs we cannot relay  on the usual variation of constants formula, as for example $w_g(\cdot)= E*\left (AB g\right )$ since   $\hat w(\zl ;w_0) $ has the decaying factor $\hat K(\zl)$  which is not present in the expression of $\hat w_g(\zl)$.  So we use the Laplace transform formulas to represent the solutions.

Thanks to the fast decaying properties of $\hat g(\zl)$ the solution $ w_g(t)$ is continuous and it is given by is 
\begin{equation}\ZLA{formGWperApproCOntro}
w_g (t )=\frac{1}{2\pi}\int _{-\ZIN}^{+\ZIN} e^{i\ZOM t} \left [ i\ZOM \hat K(i\ZOM)I-J(i\ZOM) A\right ]^{-1}
AB\left [J(i\ZOM)\hat g(i\ZOM)\right ]\ZD\ZOM\,.
\end{equation}

  Using the fact that the reachable set is increasing, $\xi\perp \mathcal{R}(T)$ (with $T\leq+\ZIN$) if and only if the following holds:
\begin{equation}\ZLA{eq:condIORTOGsu0T}
\int _{-\ZIN}^{+\ZIN} e^{i\ZOM t}\langle  \left [ i\ZOM \hat K(i\ZOM)I-J(i\ZOM) A\right ]^{-1}
ABJ(i\ZOM)\hat g(i\ZOM),\xi \rangle\ZD\ZOM=0\qquad \forall t\in [0,T]\,.
\end{equation}
 The     orthogonality condition can be written as
\begin{equation}\ZLA{eq:OrtoTfiniFreq}
\int _{-\ZIN}^{+\ZIN} e^{i\ZOM t}\langle \left ( i\ZOM+1 \right ) \hat g(i\ZOM)  ,      \frac{J(-i\ZOM)}{1-i\ZOM}B^* A \left [ -i\ZOM \hat K(-i\ZOM)I-J(-i\ZOM) A\right ]^{-1}   \xi  
   \rangle\ZD\ZOM=0 \,.
\end{equation}
Both the sides of the crochet belong  to $H^2(\Pi_+)$ (see the Appendix~\ref{appe:INFORM} for the definition and key properties of the Hardy spaces $H^2(\Pi_+)$. The minus sign in front of $i\ZOM$ in the second factor of the inner product is explained there).

The reason for inserting the factor $1/(1-i\ZOM)$ is that in this way also the second factor of the inner product is integrable on the imaginary axis, thanks to the inequalities~(\ref{eq:ineqMODIFICATE}), and so   it is the Laplace transformation of a continuous function, while  $(1+i\ZOM)\hat g(i\ZOM)$ still decays faster then any $1/|\ZOM|^k$ for every $k$.

Note that  $(1+\zl)\hat g(\zl)$ is the Laplace transform of $g+g'$ since $g(0)=0$.

Let us introduce $L(t)\xi$, the inverse Laplace transform of 
\[
 \frac{J(\zl)}{1+\zl}B^*A \left [ \zl \hat K(\zl)I-J(\zl) A\right ]^{-1}   \xi\,.
\]
The orthogonality condition~(\ref{eq:OrtoTfiniFreq}) at time $T$ takes the form (see the Appendix~\ref{appe:INFORM})
\begin{equation}
\ZLA{eq:OrtoTfiniDOMItempo}
\intt \ZL g(s)+g'(s), L(t-s)\xi\ZR\ZD\zl=0\qquad 0\leq t\leq T\leq +\ZIN\,.
\end{equation}

In order to examine the ortogonality condition, we shall use the following lemma, whose proof is at the end of this section:
\begin{Lemma}\ZLA{Lemma:prelimPERcontRappro}
Let either $U= L^2(\ZOMq_a)$ and $U_0=\mathcal{D}(\ZOMq_a)$ or $U= L^2(\Gamma_a)$ and  $U_0=\mathcal{D}(\Gamma_a)$. Let
\[
H_T=\left \{ g+g'\,,\ g\in \mathcal{D}(0,T;U_0)\right \}\,.
\]
We have:
\begin{enumerate} 
\item\ZLA{Lemma:prelimPERcontRapproDUE} if $T<+\ZIN$ then $H_T^{\perp}=\left \{ e^t u_0\,,\ u_0\in U\right \}$;
\item\ZLA{Lemma:prelimPERcontRapproUNO} If $T=+\ZIN$ then $H_{\ZIN}$ is dense in $L^2(0,+\ZIN;U)$ and so $\mathcal{H}_{\ZIN}=\left \{(1+\zl)\hat g(\zl)\,,\ g\in H_{\ZIN}\right \}$ is dense in $H^2(\Pi_+;U)$.
\end{enumerate}
\end{Lemma}

\emph{Now we characterize $\mathcal{R}_{\ZIN}^\perp$.} The proof relays on the second statement in lemma~\ref{Lemma:prelimPERcontRappro}.  
\begin{Lemma}\ZLA{Lemma:PrimoPassoproofCaraOrtinfi}
We have $\xi\perp \mathcal{R}(\ZIN)$ if and only if 
\[
\frac{J(\zl)}{\zl+1}B^* A\left (\zl\hat K(\zl)I-J(\zl) A\right )^{-1}\xi=0 \quad {\rm i.e.}\quad 
B^* A\left (\zl\hat K(\zl)I-J(\zl) A\right )^{-1}\xi=0\,. 
\]
\end{Lemma}
\zProof The two forms of the conditions are equivalent since $J(\zl)$ is not zero. 
The first form is obtained directly  from lemma~\ref{Lemma:prelimPERcontRappro} when the control acts on the boundary, since in this case $g=f$
and the orthogonality condition is that $L(t)\xi=0$ for $t>0$, i.e. $\hat L(\zl)\xi=0$ in $\zreal\zl>0$.
Analogously, the second form is immediate in the case of distributed controls, since in this case $J(\zl)\hat g(\zl)=\hat F(\zl)$.\zdia
 
\emph{This is the characterization we searched in the first step of the proof.}

\begin{Remark}\ZLA{sectApprox:RemaASSOC}
{\rm
Note that this orthogoanlity condition in the case of the associated heat equation is the well known condition
\[
B^*A\left (\zl I-A\right )^{-1}\xi=0\quad \zreal\zl>0\qquad {\rm i.e.}\qquad B^*A e^{At}\xi=0\quad t>0\,.
\]
It is known, (see~\cite{TuksnakWeiss}) that this condition implies $\xi=0$.\zdia
}
\end{Remark}

\emph{Now we proceed with the second step of the proof,} i.e. we prove that $  \mathcal{R}(\infty) ^\perp =0$. First we prove:

\begin{Lemma}\ZLA{lemmaAPProCONTRO}
The following holds:
\begin{enumerate}
\item\ZLA{item1:lemmaAPProCONTRO}
if $\xi\perp \mathcal{R}(\infty)$ then $B^*\xi=0$.
\item\ZLA{item2:lemmaAPProCONTRO} the linear space $\mathcal{R}(\ZIN)^\perp$ is invariant under $A^{-1}$: $A^{-1}\mathcal{R}(\ZIN)^\perp\subseteq \mathcal{R}(\ZIN)^\perp $ 
\end{enumerate}
\end{Lemma}
\zProof
From lemma~\ref{Lemma:PrimoPassoproofCaraOrtinfi} we have that
 $\xi\perp \mathcal{R}(\ZIN)$ if and only if
\begin{multline*}
 J(\zl)B^*A\left (\zl\hat K(\zl) I-J(\zl) A\right )^{-1}\xi\\=-B^*\xi+\zl\hat K(\zl)\left [ B^*\left (\zl \hat K(\zl) I-J(\zl)A\right )^{-1}\xi\right ]=0\,.
\end{multline*}

Now we compute the limit for $\zl\to 0$. Inequalities~(\ref{eq:ineqMODIFICATE}) show that the bracket is bounded and item~\ref{itemAssuKERNEL2ADDE2} in the set of the assumptions~\ref{itemAssu} is that the factor $\zl\hat K(\zl)$ tends to zero. Hence $B^*\xi=0$.

Once $B^*\xi=0$ is known we have also
\[
0=B^*\left (\zl \hat K(\zl) I-J(\zl)A\right )^{-1}\xi=B^*A\left (\zl \hat K(\zl) I-J(\zl)A\right )^{-1}\left [A^{-1}\xi\right ]
\]
and so $A^{-1}\xi$ is orthogonal to $\mathcal{R}(\ZIN)$ when $\xi\perp \mathcal{R}(\ZIN)$, i.e. $\left [\mathcal{R}(\ZIN)\right ]^\perp 
$ is invariant under $A^{-1}$.\zdia

 {Note that in the proof of the second statement of lemma~\ref{lemmaAPProCONTRO} we used $\zl\hat K(\zl)\not\equiv 0$, which follows from item~\ref{itemAssuKERNEL2ADDE1} in the set of the assumptions~\ref{itemAssu}.}

 \emph{Now we complete the second step of the proof}  i.e. we prove that controllability of the associated heat equation implies that if $\xi\perp \mathcal{R}(\ZIN)$ then $\xi=0$.

The operator $A^{-1} $ is selfadjoint bounded so that 
if  $\left [\mathcal{R}(\ZIN)\right ]^\perp \neq 0$ then 
$A^{-1}$ must have an eigenvalue in this invariant space:
\[
A^{-1}\xi_0=\frac{1}{\mu_0}\xi_0\quad {\rm i.e.}\quad A\xi_0=\mu_0\xi_0\qquad \xi_0\in \left [\mathcal{R}(\ZIN)\right ]^\perp \quad \mbox{and so $B^*\xi_0=0$}
\]
(note that the eigenvalue is not zero since $A^{-1} $ is invertible).  
Then we have
\[
B^*A\left [\zl I-A\right ]^{-1}\xi_0=\mu_0\frac{1}{\zl-\mu_0}B^*\xi_0=0
\]
and this is the condition that $\xi_0$ is orthogonal to the reachable set  of the associated heat equation, see Remark~\ref{sectApprox:RemaASSOC}. Hence, it must be $\xi_0=0$ i.e. \emph{we have proved the statement in item~2: $\mathcal{R}(\ZIN)$ is dense in $L^2(\ZOMq)$.}

Finally, \emph{we pass to the third step of the proof.}  We need an explicit expression of the operator $L(t)$ in~(\ref{eq:OrtoTfiniDOMItempo}). This is simply derived from the equality
\begin{multline*}
\frac{1}{1+\zl}B^*J(\zl)A\left (\zl\hat K(\zl)I-J(\zl)A\right )^{-1}\xi=-\frac{1}{\zl+1}B^*\xi\\
+B^*\left (\zl\hat K(\zl)I-J(\zl)A\right )^{-1}\hat K(\zl)\xi-\frac{1}{\zl+1}B^*\left (\zl\hat K(\zl)I-J(\zl)A\right )^{-1}\hat K(\zl)\xi
\end{multline*}
which shows
\[
L(t)\xi=-e^{-t}B^*\xi+B^*w(t;\xi)-B^*\intt e^{-(t-s)} w(s;\xi)\ZD s
\]
(we recall that $w(\cdot;\xi)$ is the solutions when the controls are put equal zero).

Now we treat separately the case of the boundary and distributed controls. In the case of the  boundary controls, $g=f$ and we use directly
statement~\ref{Lemma:prelimPERcontRapproDUE} in lemma~\ref{Lemma:prelimPERcontRappro} to deduce that   $\xi\perp \mathcal{R}(T)$ if and only if there exists $u_0$ such that
\begin{equation}\ZLA{eq:IntePerOrtTApproConttT}
-e^{-t}B^*\xi+B^*w(t;\xi)-B^*\intt e^{-(t-s)} w(s;\xi)\ZD s=e^{t}u_0\,.
\end{equation}
In fact, the orthogonality condition~(\ref{eq:OrtoTfiniDOMItempo}) is the condition that $L(\cdot)\xi\perp H_T$ where $H_T$ is defined in lemma~\ref{Lemma:prelimPERcontRappro}.

We compute~(\ref{eq:IntePerOrtTApproConttT}) for $t=0$ and we see that it must be $u_0=0$, i.e. in the case of boundary controls
 we have $\xi\perp \mathcal{R}_T$ if and only if 
\begin{equation}\ZLA{eq:VoltePERortog}
 B^*w(t;\xi)-B^*\intt e^{-(t-s)} w(s;\xi)\ZD s=e^{-t}B^*\xi\quad t\in[0,T]\,.
\end{equation}
 
The solution of the Volterra integral equation 
\[
y(t)-\intt e^{-(t-s)}y(s)\ZD s=e^{-t}B^*\xi
\]
is constant, equal to $B^*\xi$ so that
\begin{equation}\ZLA{eq:SecAPPROcondiBstWeqCONST}
B^*w(t;\xi)=B^*\xi\,,\qquad t\in [0,T]\,.
\end{equation}
 Statement~\ref{teo:EsiEVOLOPERitem3} in theorem~\ref{teo:EsiEVOLOPER} and the fact that $B$ is a bounded operator  show  that 
 $B^*w(t;\xi)$ admits a holomorphic extension to a sector surrounding the positive real axis, and so if it is constant on $[0,T]$ it is constant for every $t>0$. We compute the Laplace transformation of both the sides of the equality~(\ref{eq:SecAPPROcondiBstWeqCONST}) on $[0,+\ZIN)$ and we find
 \[
\zl\left [ B^*\left (\zl\hat K(\zl) I-J(\zl)A\right )^{-1} \hat K(\zl)\xi\right ]=B^*\xi\,. 
 \]
 We use again that $\lim _{\zl\to 0^+}\zl\hat K(\zl)=0$ and we see that $B^*\xi=0$.
 
 Then, if $\xi\perp \mathcal{R}_T$ we have also
 \begin{multline*}
0=B^*\left (\zl\hat K(\zl)\right )\left (\zl \hat K(\zl) I-J(\zl) A\right )^{-1}\xi\\
=B^*\xi+B^*J(\zl) A\left (\zl \hat K(\zl) I-J(\zl) A\right )^{-1}\xi
 \end{multline*}
 and so also
 \[
B^*A\left (\zl \hat K(\zl) I-J(\zl) A\right )^{-1}\xi\,. 
 \]
 We compare with lemma~\ref{Lemma:PrimoPassoproofCaraOrtinfi} and we see that if $\xi\perp \mathcal{R}_T$ then $ \xi\perp \mathcal{R}_\ZIN$ and so $\xi=0$ from the second step of the proof.
 
 Now we prove that if $\xi\perp\mathcal{R}(T)$ then condition~(\ref{eq:VoltePERortog}) holds also in the case of distributed controls. In fact in this case $g=(\ZDE-R)*F$ where $\ZDE$ is   the Dirac  delta and $g'=(\ZDE-R)*F'$ since $F(0)=0$. So, the orthogonality condition~(\ref{eq:OrtoTfiniDOMItempo}) is

 \[
-(\ZDE-R)*\chi B^*\xi+y-\chi*y=e^{t} u_0 
 \]
 where
 \[
\chi(t)=e^{-t}\,,\qquad y=(\ZDE-R)*B^*w(\cdot;\xi) \,.
 \]
 Computing with $t=0$ we get again $u_0=0$ and the orthogonality condition takes the form
  \[
(\ZDE-R)  \left [-\chi B^*\xi+B^*w -\chi* B^*w \right ]=0\,,\qquad w=w(\cdot;\xi)\,.
  \]
 We compute the convolution with $\ZDE+N$ and we see that $B^*w(t;\xi)$ solves~(\ref{eq:VoltePERortog}). Now we proceed as in the case of boundary controls to conclude that $\xi=0$.
 
In order to complete the proof  of  the statement in item~\ref{item:APPRX:TeoINTRIUZcontroll} of theorem~\ref{TeoINTRIUZcontroll} we prove lemma~\ref{Lemma:prelimPERcontRappro}.

 \paragraph{The proof of lemma~\ref{Lemma:prelimPERcontRappro}}
   
We prove the lemma in the case $U=L^2(\ZOMq_a)$. The proof in the case $U=L^2(\Gamma_a)$ is similar.

First we recall the following facts.
Let $f\in L^2(0,T;U)$. By definition, $f\in W^{1,2}(0,T;U)$ when we have also $f'\in W^{1,2}(0,T;U)$. The derivative is defined in the sense of the distributions, i.e.
\begin{equation}\ZLA{eq:SecAppro:ProvaLemmaOrTo}
\intT \langle f(t),\phi'(t)\rangle_U \ZD t=-\intT \ZL f'(t),\phi(t)\ZR_U \ZD t 
\end{equation} 
for \emph{any} function $\phi\in \mathcal{D}(0,T;U)$.  
In the special case under consideration, \emph{any} function $\phi\in \mathcal{D}\left (0,T;L^2(\ZOMq_a)\right )$. Note that if  $\phi\in \mathcal{D}\left (0,T;L^2(\ZOMq_a)\right )$  then we have also $\phi\in W^{1,2}\left (0,T;L^2(\ZOMq_a)\right )$.

Every $\phi\in W^{1,2}\left (0,T;L^2(\ZOMq_a)\right )$ can be approximated (in the norm of $W^{1,2} $) by a sequence $\{\phi_n\}\in \mathcal{D}\left (\ZOMq_a\times(0,T)\right )$. Hence, 
$f\in W^{1,2}(0,T;U)$ when equality~(\ref{eq:SecAppro:ProvaLemmaOrTo}) holds for every $\phi\in  \mathcal{D}\left (\ZOMq_a\times(0,T)\right )$.

 Now we prove the lemma. 
 The condition $f\perp H_T$ is the condition
\[
\intT\ZL f(t),g(t)\ZR_{L^2(\ZOMq_a)} \ZD t=-\intT\ZL f(t),g'(t)\ZR_{L^2(\ZOMq_a)} \ZD t
\]
 for every $g\in \mathcal{D}\left (\ZOMq_a\times (0,T)\right )$.
 
 It follows that
 \[
f'=f\quad \mbox{so that}\quad f(t)=e^{t} f(0)\,. 
 \]
 This proves the first statement.
 
 In order to prove the second statement we note that if $f\perp H_{\ZIN} $ then $f\perp H_T$ for every $T>0$ and so  
 both the following properties must hold: 
  $f(t)=e^{t} f(0)$ and also 
  $ f\in L^2\left (0,+\ZIN;L^2(\ZOMq)\right )$. 
  So, it must be $f(0)=0$, hence $f=0 $.\zdia

\section{\ZLA{sec:LacCONTRoZERO}Controlability to the target $0$ and reducibility to the standard heat equation}

In this section we are going to prove the existence of initial conditions $w_0$ which cannot be steered to hit the target $0$ at any time $T\notin\mathcal{Z}\subseteq (0,+\ZIN)$, an exceptional set of times which does not have accumulation points in $(0,+\ZIN)$, unless system\rmodel\ can be reduced to the standard heat equation  in the sense specified below.  

The initial condition $w_0$ can be controlled to hit the target $0$ when  
\begin{equation}
\ZLA{VhaNOcontZEROini}
\begin{array}{l}
\mbox{distributed control:  $0\in \mathcal{R}_d(w_0;T) $ i.e.  $ E(T)w_0\in   \mathcal{R}_d(0,T)$}\\
  \mbox{boundary control: $0\in \mathcal{R}_b(w_0;T)$ i.e. $
E(T)w_0\in  \mathcal{R}_b(0,T)$}\,.
\end{array}
\end{equation}
The definition of the set $\mathcal{Z}$ is suggested by the following representation of $E(t)w_0$:
  Using $I=A^{-1}A$ we see that, for $T>0$,
\begin{multline}\ZLA{eq:elabINIcondPerContzero}
2\pi i E(T)w_0\\
=A^{-1}\int_{G_\ZEP} e^{\zl T} \frac{\hat K(\zl)}{J(\zl)} \left [ J(\zl)A-\zl \hat K(\zl)I+\zl \hat K(\zl)I\right ]\left (\zl \hat K(\zl) I-J(\zl) A\right )^{-1}w_0\ZD\zl\\
 =-A^{-1}\underbrace{
 \left [\int _{G_\ZEP} e^{\zl T}\frac{\hat K(\zl)}{J(\zl)} \ZD\zl \right ]
 }_{ \mbox{we want} \ \neq 0  }w_0\\
 +\underbrace{A^{-1}\left [\int _{G_\ZEP}e^{\zl T} \frac{\zl \hat K^2(\zl)}{J(\zl)} \left (\zl \hat K(\zl) I-J(\zl) A\right )^{-1}w_0 \ZD\zl\right ] }_{\in \Dom \, A^2}\,.
\end{multline}
The second bracket is an element of $\Dom\, A$ thanks to the assumed asymptotic properties of $\hat K(\zl)$ and $\hat N(\zl)$, item~\ref{itemAssuKERNEL2ADDE} in the set of the assumptions~\ref{itemAssu}, \emph{and so $E(T) w_0\notin \Dom\,A^2$ if the first bracket is nonzero and $w_0\notin\Dom\,A$.}

So we define
\[
\mathcal{Z}=\left \{ t\,:\quad \int _{G_\ZEP} e^{\zl t}\frac{\hat K(\zl)}{J(\zl)} \ZD\zl=0\ \right\}\,.
\]
 
We are going to see  that $\mathcal{Z}$ is the set  of the real zeros of a function which is holomorphic in a sector surrounding the real positive axis.   So, it does not have accumulation points in $t>0$, unless the function is identically zero. In this case we shall see that  system~(\ref{eq:modello})  can be reduced to the associated heat equation. In fact:
\begin{Lemma}\ZLA{Lemma:NOaccumPOINT}
The function 
\[
  \Psi(z)=\frac{1}{2\pi i} \int _{G_\ZEP} e^{\zl z} \frac{\hat K(\zl)}{J(\zl)} \ZD\zl
\]
is holomorphic in the sector $\Sigma_{\theta+\pi/2}$ ($\theta$ is as in Assumption~\ref{itemAssu} item~\ref{itemAssuKERNEL2}) and it is zero if and only if $\frac{\hat K(\zl)}{J(\zl)}$ is a   constant.
If $\frac{\hat K(\zl)}{J(\zl)}\equiv c$ then $c>0$.
\end{Lemma}
\zProof
Note that $ {\hat K(\zl)}/{J(\zl)}$ is bounded in the set 
$\Sigma_{\zthe+\pi/2}\cap\{\zl\,:\ |\zl|>a\}$ (any $a>0$) so that the improper integral converges uniformly on compact sets of $\Sigma_{\theta+\pi/2}$, i.e. $\Psi(z)$ is holomorphic on $\Sigma_{\theta+\pi/2}$.

If $ {\hat K(\zl)}/{J(\zl)}\equiv c$ then the integral is zero, as seen by integrating on the closed circuit $G_\ZEP\cap \cal\{\zl\,:\  |\zl|<R\}$, closed by an arc of radius $R$ in the left half plane, and passing to the limit for $R\to+\ZIN$.

Conversely, we prove that if $\Psi(z)\equiv 0$ then $ {\hat K(\zl)}/{J(\zl)}$ is constant and in fact \emph{it is a positive constant}. 

As in~\cite[p.~59]{PruessLIBRO1993}, we introduce
\[
\Phi(z)= \frac{1}{2\pi i} \int _{G_\ZEP} e^{\zl z} \frac{1}{\zl^2}\frac{\hat K(\zl)}{J(\zl)} \ZD\zl\,.
\]
Also this function is holomorphic on $\Sigma_{\theta+\pi/2}$.

Thanks to the fast decaying property of the exponential,

\begin{equation}\ZLA{eq:PhiComEpOlInOmIo}
\Phi''(z)=\Psi(z)\equiv 0\quad \mbox{i.e.}\quad \Phi(z)=\Phi_0+\Phi_1 z  \,. 
\end{equation}
 So we have (see formula~(\ref{Eq:Appe:AtitTraslata}))
\begin{multline*}
\Phi(t)=\Phi_0+\Phi_1 t=\frac{1}{2\pi i} \int _{G_\ZEP} e^{\zl t} \frac{1}{\zl^2}\frac{\hat K(\zl)}{J(\zl)} \ZD\zl\\\
=
\frac{1}{2\pi }\int _{c- i\ZIN}^{c+ i\ZIN} e^{ (c+ i\ZOM) t} \frac{\hat K(c+ i\ZOM)}{(c+i\ZOM)^2J(c+ i\ZOM)}\ZD\ZOM\,,\qquad c>\ZEP\,.
\end{multline*}
The last equality is justified by the facts that, from the properties~\ref{itemAssuKERNEL2ADDE0}
and~\ref{itemAssuKERNEL3} in the set   of the assumptions~\ref{itemAssu},   $\hat K(\zl)$ is bounded on $G_\ZEP$ and in the region to the right of it and $J(\zl)\neq 0$.
It follows that the Laplace transform of $\Phi(t)$ is 
\[
\hat \Phi(\zl) =\frac{\hat K(\zl)}{\zl^2J(\zl)}\qquad \zreal\zl>0.
\]
Using the expression of $\Phi(t)$ in~(\ref{eq:PhiComEpOlInOmIo})  we have
\[
\frac{\hat K(\zl)}{\zl^2J(\zl)}=\Phi_0\frac{1}{\zl}+\Phi_1\frac{1}{\zl^2}\quad \mbox{and so $\frac{\zl \hat K(\zl)}{J(\zl)}=\Phi_0\zl^2+\Phi_1\zl\,.$}
\]
Now we invoke the assumption~\ref{itemAssu} item~\ref{itemAssuKERNEL4}: the map
$\zl
 \mapsto \zl\hat K(\zl)/J(\zl) 
 $
 transforms $\Sigma_{\zthe+\pi/2} $ to $\Sigma_{\zthe_A+\pi/2}$. This is   possible only if $\Phi_0=0$
 and $\Phi_1\geq 0$ (and in fact $\Phi_1>0$ since $K\neq 0$)   as wanted.\zdia

And so:
\begin{Lemma}\ZLA{Lemma:ReducASSOCheatEQ}
The set $\mathcal{Z}$ has accumulation points in $t>0$, and so  $\mathcal{Z}=(0,+\ZIN)$, if and only if Eq.~(\ref{eq:modello}) \emph{is the associated heat equation.}
\end{Lemma}
\zProof  
We proved that if   $\mathcal{Z}=(0,+\ZIN)$ then
\[
\hat K(\zl)/J(\zl)\equiv \Phi_1>0
\quad{\rm i.e.}\quad  \hat K(\zl)=\Phi_1\left (1+\hat N(\zl)\right )\,,\qquad \Phi_1>0\,.
\]

Then, Eq.~(\ref{eq:modello})
takes the form
\[
\left (w' -\frac{1}{\Phi_1}\Delta w\right )+\intt N(t-s)\left (w' (s)-\frac{1}{\Phi_1}\Delta w(s)\right )\ZD s=\frac{1}{\Phi_1} F\,.
\]
Let $R(t)$ be the resolvent kernel of $N(t)$. Then we have
 \[
w'(t) -\frac{1}{\Phi_1}\Delta w(t)= \frac{1}{\Phi_1} F-\intt R(t-s) \frac{1}{\Phi_1} F(s)\ZD s \,,\qquad \Phi_1>0
 \]
 (with the same boundary control $f$). This is the associated heat equation (with a diffusion coefficient $ 1/\Phi_1 $ possibly not equal $1$, but it is clear that this has no importance. If we want also the   coefficient to be equal $1$ this can be achieved by changing the time scale).\zdia
 
   Lemma~\ref{Lemma:ReducASSOCheatEQ} suggests the following definition: 
 
 \begin{Definition}\ZLA{defin:ReducHEAt}   {\rm
 System\rmodel\ is reducible to the standard heat equation when there exists $c >0$ such that $K(t)=c \left (\ZDE+N(t)\right )$ where $\ZDE$ is the  Dirac  delta.
 }
 \end{Definition}

Now we prove:

\begin{Theorem}\ZLA{Teo:LackZEROconTRO}
Let $\mathcal{Z}  \neq (0,+\ZIN)$ (i.e.   the system\rmodel\ is not reducible to the standard heat equation) and let
  $T\notin\mathcal{Z}$. Then:
\begin{enumerate}
 \item\ZLA{item2:Teo:LackZEROconTRO}  
if $\ZOMq\setminus {\rm cl}\ZOMq_a\neq\emptyset$ then there exist initial conditions $w_0$ (which do not depend on $T\notin\mathcal{Z}$) such that $E(T)w_0\notin \mathcal{R}_d( 0;T)$;
\item\ZLA{item1:Teo:LackZEROconTRO} there exist initial conditions $w_0$ (which do not depend on $T\notin\mathcal{Z}$) such that $E(T)w_0\notin \mathcal{R}_b( 0;T)$.
\end{enumerate}
\end{Theorem} 
\zProof 
The proof in the case of distributed controls is immediate:  if $F\in \mathcal{D}\left ( \ZOMq_a\times(0,T)\right )$ then  
$w(T;F)\in \Dom\,A^k$ for every $k$,   see the statement in item~\ref{teo:EsiSOLUaffine:itemCONdistri} of theorem~\ref{teo:EsiSOLUaffine}. Hence,  the inclusion $w(T;w_0)\in \mathcal{R}_d(0,T)$ does not hold for example when $w_0\notin \Dom\, A$, thanks to the representation~(\ref{eq:elabINIcondPerContzero}) of $E(T)w_0$.

The proof of lack of controllability under boundary controls follows similar ideas. Thanks to the regularity of $f(t)$, we have that $(\zl+1)^k\hat f(\zl)$ is bounded in the right half plane 
for every $k\geq 0$ and formula~(\ref{eq:TrasfoRISOLTA}) gives the following expression for $w_f(t)$ (here $R_c$ is the vertical line $\zreal\zl =c$ for an arbitrary fixed $c>0$)
\begin{multline*}
w_f(t)=Gf(t)\\
-\frac{1}{2\pi i}\underbrace{\int _{R_c} e^{\zl t} \frac{\zl\hat K(\zl)}{(\zl+1)^k}\left (\zl\hat K(\zl)I-J(\zl)A\right )^{-1} G\left [(\zl+1)^k\hat f(\zl)\right ]\ZD\zl}_{\tilde w(t)}
\end{multline*}
where the path of integration   $R_c$ is the vertical line $\zreal\zl =c$ for an arbitrary fixed $c>0$,
\[
R_c=c+i\ZOM\,,\qquad \ZOM\in\zzr\,.
\]
We shall see that $\tilde w(t)$ takes values in $\Dom(-A)^{2-\ZSI_0}$ for every $t$, in particular for $t=T$.  Accepting  this fact and using $f(T)=0$ we then have $w_f(T)\in \Dom(-A)^{2-\ZSI_0}$ and the inclusion $E(T)w_0\in \mathcal{R}_b(T)$ cannot hold when $w_0\notin \Dom(-A)^{1-\ZSI_0}$.

In order to complete our argument we  prove the stated regularity of $\tilde w(t)$. This follows  from the following equality:
\begin{multline*}
A\tilde w(t)\\=  \int _{R_c}e^{\zl t}\frac{\zl }{(\zl+1)^k}
 \frac{\hat K(\zl)}{J(\zl)}
 \left \{
-I+\zl\hat K(\zl)\left (\zl\hat K(\zl) I-J(\zl)A\right )^{-1}
\right \} G\left [(\zl+1)^k \hat f(\zl)\right ]\ZD\zl\\
=- \underbrace{G\int_{R_c}e^{\zl t}
\frac{\zl\hat K(\zl)}{J(\zl)} 
\hat f(\zl)\ZD\zl}_{\in \Dom(-A)^{\ZSI_0}
 }\\  +  
  \underbrace{
  \int _{R_c}e^{i\zl t} \left \{\frac{\zl^2\hat K(\zl)^2}{(\zl+1)^k}
  \left (\zl\hat K(\zl) I-J(\zl) A\right )^{-1}\right \}G\left [(\zl+1)^k\hat f(\zl)\right ]\ZD\zl
 }_{\in\Dom A}\,.
\end{multline*}

The proof of the theorem is now complete.\zdia

 \begin{Remark}
 {\rm
 Note that the result   concerning  distributed controls is stronger than stated  since  we did not use the full strength of the assumption that the control is localized in $\ZOMq_a$ (which is the case of interest in practice) but we used solely the fact that it is possible to control $w_0$ to hit the target zero with distributed controls which are smooth enough. Also the full strength of the condition that the controls are of class $C^\ZIN$ has not been used.\zdia
 }
 \end{Remark}

 \appendix
 \section{Appendices}
 \subsection{\ZLA{appe:INFORM}Information on the Laplace transform and Hardy spaces}
Let $X$ be any separable Hilbert space. We use $\Pi_+$ to denote the right half plane, $\Pi_+=\{\zreal\zl>0\}$. The Hardy space $H^2(\Pi_+;X)$ is the space of the   functions $F(\zl)$ which are holomorphic in $\zreal\zl>0$ and such that
\[
\|F\|^2_{H^2(\Pi_+;X)}=\sup _{x>0}\int_{-\ZIN}^{+\ZIN} \|F(x+iy)\|_{X}^2\ZD y<+\ZIN\,.
\]
The   linear space $H^2(\Pi_+;X)$ endowed with this norm is a Hilbert space and the Laplace transformation $f\mapsto \hat f=F$ is bounded and boundedly invertible between $L^2(0,+\ZIN;X)$ and $H^2(\Pi_+;X)$.

We need the following additional pieces of information:
\begin{itemize}
\item if $F\in H^2(\Pi_+;X)$ then $\lim _{x\to 0^+}F(x+iy)=F(iy)$ exists a.e. and exists in the sense of $L^2(i\zzr)$, i.e. $\lim _{x\to 0^+}\|F(x+iy)-F(iy)\|_{L^2(-\ZIN,+\ZIN;X)}= 0$.  
\item
The space $\{ F(iy)\,,\ F(x+iy)\in H^2(\Pi_+;X)\}$ 
is a closed subspace of $L^2(-\ZIN,+\ZIN;X)$ and it turns out that 
\begin{multline*}
\langle F,G\rangle _{H^2(\Pi_+;X)}=\sup _{x>0}\left \{\int _{-\ZIN}^{+\ZIN} \bar G(x+iy)F(x+iy)\ZD y\right \}
\\
=\int  _{-\ZIN}^{+\ZIN} \bar G( iy)F( iy)\ZD y
\end{multline*}
(i.e. $H^2(\Pi_+;X)$ and such closed subspace of  $L^2(-\ZIN,+\ZIN;X)$ are isometric).

\item In the special case of a real valued function $f(t)$, $\overline{\hat f}(x+iy)=\hat f(x-iy)$ and so, if $f$ and $g$ are real valued and $F=\hat f$, $G=\hat g$, we have
 
\begin{multline*}
\langle F,G\rangle _{H^2(\Pi_+;X)}=\sup _{x>0}\left \{\int _{-\ZIN}^{+\ZIN} \bar G(x+iy)F(x+iy)\ZD y\right \}
\\=\int  _{-\ZIN}^{+\ZIN}  G(- iy)F( iy)\ZD y\,.
\end{multline*}
 \item Let $f$ and $g$ belong to $L^2(0,+\ZIN;X)$. Then,
 \[
h(t)=\intt \ZL f(t-s),g(s)\ZR \ZD s\in L^\ZIN(0,+\ZIN)\cap C([0,+\ZIN)) 
 \] 
and we have:
 \begin{itemize}
 \item
 The Laplace transform of $h(t)$ is 
 \[
\hat h(\zl)=\ZL \hat f(\zl),\overline{\hat g}(\zl)\ZR_X\,.
 \]
 
 \item In the special case that $f(t)$ and $g(t)$ take values in a \emph{real} Hilbert space $X$ then $h(t)$ is real valued   and
 \[
\hat h(i\ZOM)= \ZL \hat f(i\ZOM),\overline{\hat g}(i\ZOM)\ZR_X=\ZL \hat f(i\ZOM),  \hat g(-i\ZOM)\ZR_X\,.
 \]
 \end{itemize}
\item when $f\in L^2(0,+\ZIN;X)$ i.e. $\hat f\in H^2(\Pi_+;X)$, the usual formula of the inverse Laplace transform holds
\begin{equation}\ZLA{eq:FormAtiTraLaplaImagAX}
f(t)=\frac{1}{2\pi} \int _{-\ZIN}^{+\ZIN}e^{i\ZOM t}\hat f(i\ZOM)\ZD\ZOM
\end{equation}
in the sense that the integral computed on $[-T,T]$ converges to $f$ in $ L^2(0,+\ZIN;X)$ when $T\to +\ZIN$. If it happens that $\hat f(i\ZOM)\in L^1(-\ZIN,+\ZIN;X)$ then $f(t)$ is continuous and the convergence (for $T\to+\ZIN$) is uniform on compact subsets of $[0,+\ZIN)$.
\item If there exists $\zaa\in\zzr$ such that $F(\zl+\zaa)\in H^2(\Pi_+;X)$ then $F(\zl)=\hat f(\zl)$ where the function $f(t)$ is given by
\begin{equation}\ZLA{Eq:Appe:AtitTraslata}
f(t)=\frac{1}{2\pi}\int _{c-i\ZIN}^{c+i\ZIN} e^{(c+i\ZOM) t} F(c+i\ZOM)\ZD\ZOM=
\frac{1}{2\pi i}\int _{R_c}ve^{\zl  t} F(\zl)\ZD\zl
\end{equation}
where $R_c=c+i\ZOM$, $-\ZIN<\ZOM<+\ZIN$.
The number $c$ is any real number larger then $a$.
 \item When $f\in\mathcal{D}\left ( (0,+\ZIN)\times X\right )$ 
  its Laplace transform $\hat f(\zl)$ decays faster then $1/|\zl|^k$ when $|\zl|\to+\ZIN$ in $\zreal\zl\geq 0$, for every $k>0$.
\end{itemize}

\subsection{\ZLA{appe:ProofsEXI}The proofs of theorems~\ref{teo:EsiEVOLOPER}   and~\ref{teo:EsiSOLUaffine}}
The proofs follows estabilished routes, see for examples~\cite{PruessLIBRO1993}. In fact, similar ideas are used in the study of holomorphic semigroups, see~\cite{Pazy}. 

We sketch the proofs in order to see the role of the assumptions, in particular of the assumptions on $\hat K(\zl)$.

\begin{figure}[h]
\caption{\label{fig:path}The path  of integration $G_\ZEP$. The dotted angle is $\Sigma_{\zthe}$. } 
\vspace{0cm}
\begin{center}
\includegraphics[width=8cm]{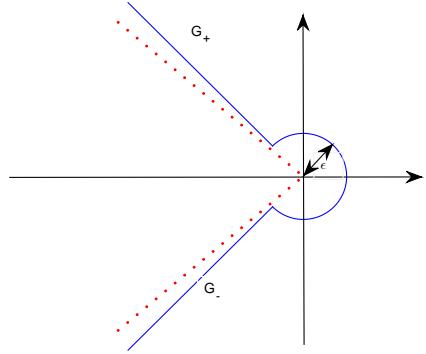}
\end{center}
\smallskip
\center{}
\end{figure}

The idea is to recover the candidate solution  $w(t;w_0)=E(t)w_0$
 as the inverse Laplace transform of $(\zl\hat K(\zl)I-J(\zl) A)^{-1}\hat K(\zl)$ but 
the restriction to the imaginary axis is not integrable so that in order to find the inverse Laplace transform we integrate   on the path   $G_\ZEP$ in~(\ref{eq:intepath}) and represented in Fig.~\ref{fig:path}.
 
We recall that the angle $\zaa\in (\pi/2,\zthe)$ is fixed.

Then we define
\[
E(z)=\frac{1}{	2\pi i}\int _{G_\ZEP} e^{z\zl} \hat K(\zl)\left (\zl\hat K(\zl)I-J(\zl)A\right )^{-1}\ZD \zl\,.
\]
Of course the improper integral is computed as the limit of integrals on $G_{R,\ZEP}=G_{\ZEP}\cap \{\zl\,\ |\zl|<R\}$. It is easy to see that the limit exists uniformly on compact subsets of the sector
\[
\Sigma_{\zthe}=\{z \,:\ |\arg z|<\zthe\,,\ z\neq 0\}
\]
and so $E(z)$ is a holomorphic operator valued function on $\Sigma_\zthe$.

The operator $AE(z)$   has similar properties. This is simply seen as follows: the operator valued function  $AE(z)$ is bounded on every compact subset of $\Sigma_\zthe$ as seen from
\[
  \hat K(\zl)A\left (\zl\hat K(\zl)I-J(\zl)A\right )^{-1} 
=-I+\frac{\zl\hat K(\zl)}{J(\zl)}\left (\frac{\zl\hat K(\zl)}{J(\zl)}I-A\right)^{-1} \,.
\]
The second addendum is bounded from the assumptions on the memory kernels and the first set of the inequalities~(\ref{eq:ineqMODIFICATE}) while 
   $e^{z\zl}$ is dominated by a decaying exponential if $z$ belongs to a \emph{compact} subset of $\Sigma_\zthe$ (in particular it is far from $0$). 
   
Let now $w\in L^2(\ZOMq)$ and $\tilde w\in\Dom A$. Then we have (the crochet denotes the inner product in $L^2(\ZOMq)$)
\[
\langle AE(z)w,\tilde w\rangle=\langle E(z)w,A\tilde w\rangle
\] 
so that $z\mapsto \langle AE(z)w,\tilde w\rangle$ is a holomorphic valued function for every $\tilde w\in\Dom A$ (which is dense in $L^2(\ZOMq)$). We use~\cite[Remark~1.38 p.~139]{KATOperturbTheoryBOOK} and we deduce  that $z\mapsto AE(z)$ is a holomorphic   operator
valued function on $\Sigma_\zthe$.

Now we prove that $E(t)$ is bounded for $t>0$. We proceed as in~\cite{DaPratoPADOVA1980}: we change the variable $\zl=\zeta/t$ and we see that $E(t)$ is given by
\[
E(t)=\frac{1}{2\pi i} \int_{G_{\ZEP t}} e^{\zeta}\left (1/t\right )\hat K (\zeta/t  )\left (\frac{\zeta}{t}\hat K(\zeta/t)I-J(\zeta/t)A\right )^{-1}\ZD\zeta\,.
\]
The straight lines of the paths $G_{\ZEP t} $ and $G_\ZEP$ coincide but the paths are closed by the circular arcs of radius respectively $\ZEP t$ and $\ZEP$ (with the same center zero). So, the integration paths differ by a path which does not enclose singularities of the integrand
and we have also
\[
E(t)=\frac{1}{2\pi i} \int_{G_{\ZEP  }} e^{\zeta}\left (1/t\right )\hat K (\zeta/t  )\left (\frac{\zeta}{t}\hat K(\zeta/t)I-J(\zeta/t)A\right )^{-1}\ZD\zeta\,.
\]
Boundedness  uniformly in $t >0$ follows since
\[
\left | \frac{1}{t}\hat K (\zeta/t  )  \right |\,\left \| 
\left (\frac{\zeta}{t}\hat K(\zeta/t)I-J(\zeta/t)A\right )^{-1}
\right \|\leq \frac{M}{|\zeta|}\,.
\]

Finally, we prove that
\[
\lim _{t \to 0+} E(t)w_0=w_0\quad \forall w_0\in L^2(\ZOMq)\,.
\]

As we already noted that $E(t)$ is a bounded function of $t$, we can confine ourselves to prove this property for $w_0=A^{-1} y_0$.

We use the first formula of the resolvent and we see that
\begin{multline*}
\frac{1}{2\pi i}\int _{G_\ZEP} e^{\zl t}\hat K(\zl)\left (\zl\hat K(\zl) I-J(\zl)A\right )^{-1} w_0\ZD \zl
\\
=\frac{1}{2\pi i}\int_{G_\ZEP} e^{\zl t} \frac{1}{\zl} A^{-1}y_0\ZD\zl+\frac{1}{2\pi i}\int_{G_\ZEP}e^{\zl t} \frac{1}{\zl}  \left (\frac{\zl\hat K(\zl)}{J(\zl)}I-A\right )^{-1}y_0\ZD \zl\,.
\end{multline*}
The first integral is $A^{-1}y_0$ from Cauchy integral formula and so it is sufficient to note that the limit of the second integral for $t\to 0+$ is equal zero.  To prove this fact, we again use the transformation $\zl=\zeta/t$ and write the second integral as
\[
\int _{G_\ZEP} e^{\zeta}\left (1/\zeta\right )\left (\frac{\zeta}{t} \frac{\hat K(\zeta/t)}{J(\zeta/t)}I-A\right )^{-1} y_0\ZD\zeta\,.
\]
The condition in   item~\ref{itemAssuASYMPkernelK} of the set of the assumptions~\ref{itemAssu} 
shows that
\[
\left \|\left (\frac{\zeta}{t} \frac{\hat K(\zeta/t)}{J(\zeta/t)}I-A\right )^{-1}\right |\leq M\frac{|J(\zeta/t)|}{|(\zeta/t)\hat K(\zeta/t)|}\leq \frac{M}{|\zeta/t|^{\zg_0}} 
\]
tends to zero when $t\to 0^+$. The integral tends to zero since $e^{\zeta }/\zeta $ tends to zero exponentially fast when $|\zeta|\to+\ZIN$, $\zeta\in G_\ZEP$.

Now we compute the Laplace transform   $\hat E(s)w_0$  in $\zreal s>0$  of $t\mapsto E(t)w_0$
and we prove that it is precisely $\hat K(\zl)\left (\zl\hat K(s\zl)I-J(\zl)A\right )^{-1}w_0$, as it must be if we want to choose $E(t) w_0$ as the definition of the solution $w(t;w_0)$:
\begin{multline*}
\hat E(s)w_0=\frac{1	}{2\pi i}\int_0^{+\ZIN} e^{-st}\left [\int _{G_\ZEP}
e^{t\zl} \hat K(\zl)\left (\zl\hat K(\zl)I-J(\zl)A\right )^{-1}w_0\ZD \zl\right ]\ZD t\\
=\frac{1}{2\pi i} \int _{G_\ZEP} \frac{1}{s-\zl	}  \hat K(\zl)\left (\zl\hat K(\zl)I-J(\zl)A\right )^{-1}\ZD \zl\\
=
 \hat K(s)\left (s \hat K(s)I-J(s)A\right )^{-1}w_0\,.
\end{multline*}
The last integral is computed as the limit of the integrals on $G_{R,\ZEP}=G_{\ZEP}\cap \{\zl\,\ |\zl|<R\}$  completed with the circular arch of radius $R$ and \emph{which intersect the right half plane $\zreal\zl>0$.} The equality follows since the  integrand has the sole singularity  $\zl=s$ in the interior region while
the contribution of the circular arch tends to zero since the integral decays as $1/\zl^2$ (we take into account the fact that the path is described in 
  the
\emph{negative } sense).  

This ends the proof of theorem~\ref{teo:EsiEVOLOPER} and we proved also the representation formula~(\ref{eq:rAPPreDIwinidata}).\zdia

The previous arguments show that it makes sense to use $w(t;w_0)=E(t)w_0$ has the mild solution of Eq.\rmodel (when $F=0$, $f=0$) and that the \emph{evolution operator $E(t)$} admits a holomorphic extension to a sector enclosing the axis $t>0$.

In order to prove theorem~\ref{teo:EsiSOLUaffine} we first examine $\left  (\zl \hat K(\zl) I-J(\zl) A\right )^{-1}\hat F(\zl)$ and we prove that when $F\in L^2\left (0,+\ZIN;L^2(\ZOMq)\right )$ this function is the Laplace transformation of a square integrable $L^2(\ZOMq)$-valued function. 

  By assumption,
 $\hat F(\zl)\in H^2\left (\Pi_+;L^2(\ZOMq)\right ) $ and 
\begin{equation}\ZLA{eq:addeTOappe1}
\left(\zl\hat K(\zl)-J(\zl) A\right )^{-1}\hat F(\zl)\in H^2\left (\Pi_+;L^2(\ZOMq)\right )\,.
\end{equation}
This follows from the fact that 
\[
\left\| \left  (\zl \hat K(\zl) I-J(\zl) A\right )^{-1}\right \| \leq     
 \frac{M}{\left 	|\zl\hat K(\zl)+\ZOM J(\zl)\right |}
\] 
(this is the  first inequality in~(\ref{eq:ineqMODIFICATE})). In fact we noted in Remark~\ref{rema:sulNONnullo} that the denominator is not zero, and so it is sufficient to note that it does not approach zero for $\zl\to 0$ and $|\zl|\to+\ZIN$ in $\zreal\zl>0$. This follows easily from the assumption~\ref{itemAssu} item~\ref{itemAssuASYMPkernelK}.

The fact that the multiplier $ \left  (\zl \hat K(\zl) I-J(\zl) A\right )^{-1}$ is holomorphic 
and bounded in the right half plane implies also that the transformation from $\hat F\in H^2\left (\Pi_+;L^2(\ZOMq)\right )$ to the function in~(\ref{eq:addeTOappe1}), as an element of $H^2\left (\Pi_+;L^2(\ZOMq)\right )$, is bounded and so, passing to the time domain, 
there exists a transformation $\mathcal{E}_d\in \mathcal{L}\left ( L^2\left (0,T;L^2(\ZOMq)\right ) \right )$ such that
the Laplace transform of $\left (\mathcal{E}F\right )(t)$ is $\left  (\zl \hat K(\zl) I-J(\zl) A\right )^{-1}\hat F(\zl)$, as we wanted to prove.

For every $T>0$ the transformation $\left (\mathcal{E}_{d,T}F\right )$ is just the restriction of $\mathcal{E}F$ to $(0,T)$  (in case that $F$ is defined only on $(0,T)$ we intend that it is extended with $0$ for $t>T$).

This justify taking the inverse Laplace transform of
 $
 \left  (\zl \hat K(\zl) I-J(\zl) A\right )^{-1} \hat F(\zl) \,.
$ as the solution 
$
w(t;F)
$.

Now we prove the representation formula~(\ref{rappreDISTRIBOUNDcontro}).
Using the formula for the Laplace transform of a convolution, it is sufficient to prove that the 
integral~(\ref{RappreEcorsivo}) converges in $L^1\left (0,T;\mathcal{L}(X)\right )$ (convergence in $C\left ( [a,T];\mathcal{L}(X)\right )$ for $a>0$ is clear) and then to compute its Laplace transform.
 
In order to prove convergence in $L^1\left (0,T;\mathcal{L}(X)\right )$ (the integral is in Bochner sense) we prove convergence of the integral of the norm. It is sufficent that we consider the integrals on $G_+$ (the integral on $G_-$ is treated analogously). The parametrization of $G_+$ is in~(\ref{eq:intepath}). We see that

\[
\left \|e^{\zl t}\left( \zl\hat K(\zl)-J(\zl) A\right )^{-1}\right \|\leq e^{-st|\cos\zaa|} \frac{M}{s^{\zg_0}}
\]
and we must prove that
 
\[
\lim_{R\to+\ZIN}\int_0^{+\ZIN} \int_R^{+\ZIN} e^{-stc}\frac{1}{s^{\zg_0}} \ZD s\,\ZD t=0\,,\qquad c=|\cos\zaa|\,.
\]
We replace   $stc=\nu$ in the inner integral and then $Rtc=\xi$ and we see that  the integral is equal to
\begin{multline*}
\frac{1}{cR^{\zg_0}} 
\int_0^{+\ZIN} \frac{1}{\xi^{1-\zg_0}} \int_\xi^{+\ZIN} e^{-\nu}\frac{1}{\nu^{\zg_0}} \ZD\nu\,\ZD \xi
\\=\frac{1}{cR^{\zg_0}} \int_0^{+\ZIN}e^{-\nu}\frac{1}{\nu^{\zg_0}}\int_0^{\nu}\frac{1}{\xi^{1-\zg_0}}\ZD\xi\,\ZD\nu
=\frac{1}{\zg_0cR^{\zg_0}}\to 0\,.
\end{multline*}

Using the Laplace transform of a convolution, formula~(\ref{rappreDISTRIBOUNDcontro}) is   by proving that  the Laplace transform of~(\ref{RappreEcorsivo}) is $\left (\zl\hat K(\zl)I-J(\zl)A\right )^{-1}$, similar to what we did above.
 
Thanks to the convergence of the integral we can exchange the order of integration in the computation of the Laplace transformation
\begin{multline*}
\int_0^{+\ZIN} e^{-\zeta t}\left [
\frac{1}{2\pi i}\int _{G_\ZEP} e^{\zl t} \left (\zl\hat K(\zl)-J(\zl) A\right )^{-1}\ZD\zl
\right ] \ZD t\\
=\frac{1}{2\pi i} \int _{G_\ZEP} \left [ \int _0^{+\ZIN}  e^{t( \zeta-\zl)}\ZD t\right ]\left (\zl\hat K(\zl)-J(\zl) A\right )^{-1}\ZD\zl\\
=\frac{1}{2\pi i} \int _{G_\ZEP}   \frac{1}{\zl-\zeta}\left (\zl\hat K(\zl)-J(\zl) A\right )^{-1}\ZD\zl=
\left (\zeta\hat K(\zeta)-J(\zeta) A\right )^{-1}
\end{multline*}
from Cauchy integral formula.

An analogous argument (based on the second set of inequalities in~(\ref{eq:ineqMODIFICATE})) shows that $f\mapsto w_f(t)$ ($t\in [0,T]$) is linear and continuous from $L^2\left (0,T;L^2(\Gamma_a)\right )$ 
to $L^2\left (0,T;L^2(\ZOMq)\right )$ and the representation formula~(\ref{rappreDISTRIBOUNDcontro}).

Now we observe that $w(t;F)$ and $w_f(t)$ are continuous $L^2(\ZOMq)$ valued functions when $F$, or $f$  are $C^\ZIN$ with compact support. In fact, in this case for every $k$ there exists $M_k>0$ such that, respectively,
\[
\left \|\hat F(\zl)\right \|_{L^2(\ZOMq)}\leq \frac{M_k}{|\zl|^k}\,,\quad 
\left \|\hat f(\zl)\right \|_{L^2(\Gamma)}\leq \frac{M_k}{|\zl|^k}
\]
so that
\[
\left (i\ZOM\hat K(i\ZOM)-J(i\ZOM)A\right )^{-1}\hat F(i\ZOM)\,,\qquad 
\left (i\ZOM\hat K(i\ZOM)-J(i\ZOM)A\right )^{-1}AG\hat f(i\ZOM)
\]
are integrable on the imaginary axis and the inverse Laplace transformations are continuous $L^2(\ZOMq)$-valued functions.

Finally, if $F\in\mathcal{D}\left (\ZOMq\times(0,+\ZIN)\right )$ then we have also $F\in C^{\ZIN}\left ([0,T],\Dom A^k\right )$ (and also $F^{(k)}(0)=0$) for every $k$ so that the representation formula~(\ref{rappreDISTRIBOUNDcontro}) shows that $w(t;F)\in C\left ([0,T];\Dom\, A^k\right )$ for every $k$.
This completes the proof theorem~\ref{teo:EsiSOLUaffine}.\zdia

Finally, we give an example which shows that $t\mapsto w(t;F)$ needs not be continuous. We use the semigroup property of the Riemann-Liouville integral:

\[
J^\zg\circ J^\ZSI=J^{\zg+\ZSI}.
\]
\begin{Example}\ZLA{EXE:nOcONTIN}
{\rm
We fix $\ZEP\in (0,1/4)$ and $\zg= \ZEP+1/2<1$.
We consider the equation
\[
J^\zg w'=\Delta w+F\qquad w(0)=w_0=0\,,\quad f=0\,.
\]
 We apply $J^{1-\zg}$ to both the sides and we get
the equation
\[
J^1w'=w(t)=J^{1-\zg}\Delta w+J^{1-\zg} F\,.
\]
We choose $F=F(x,t)=F_0(t) $ and we put $c_n=\int_\ZOMq\phi_n(x)\ZD x$ where $\phi_n$ is an eigenfunction of the operator $A$ such that $c_n\neq 0$ ($-\mu_n^2$ be the eigenvalue).

We project on the eigenfunction $\phi_n$. Let $w_n(t)=\int_{\ZOMq}w(x,t)\phi_n(x)\ZD x$. We get
\[
w_n(t)+\mu_n^2 J^{1-\zg} w_n = c_n J^{1-\zg} F\,.
\]
If $t\mapsto w(\cdot ,t)$ is a continuous $L^2(\ZOMq)$ valued function, then the left hand side is a continuous function of $t$. In contrast with this, the right hand side   is not continuous at $t=T$ when, for example,
\[
F(t)=\frac{1}{\left (T-s\right )^{(1/2)-\ZEP}}\in L^2(0,T)\,.
\]
In fact in this case    we have for $t<T$
\[
(J^{1-\zg} F)(t)=\frac{1}{\Gamma(\zg)}\intt \frac{1}{(t-s)^{(1/2)+\ZEP}} \frac{1}{(T-s)^{(1/2)-\ZEP}}\ZD s\,.
\]
This is a continuous function of $t$ for $t<T$ but the   limit for $t\to T^-$ is not finite. In fact, $t-s\leq T-s$ so that
\begin{align*}
\lim _{t\to T^-}\intt \frac{1}{(t-s)^{(1/2)+\ZEP}} \frac{1}{(T-s)^{(1/2)-\ZEP}}\ZD s\geq
\lim _{t\to T^-} \intt \frac{1}{T-s}\ZD s\\
= \lim _{t\to T^-}\left [\log T-\log(T-t)\right ]= +\ZIN \,.\zdiaform
\end{align*}

}
\end{Example}

 \bibliography{bibliomemoria}{ }

\begin{thebibliography}{10}

\bibitem{BarbuDIFFINTEQ2000}
V.~Barbu and M.~Iannelli.
\newblock Controllability of the heat equation with memory.
\newblock {\em Differential Integral Equations}, 13(10-12):1393--1412, 2000.

\bibitem{BaumeisterJMAA1983}
J.~Baumeister.
\newblock Boundary control of an integro-differential equation.
\newblock {\em J. Math. Anal. Appl.}, 93(2):550--570, 1983.

\bibitem{ColemanGurtinZEITSCH1967}
B.~D. Coleman and M.~E. Gurtin.
\newblock Equipresence and constitutive equations for rigid heat conductors.
\newblock {\em Z. Angew. Math. Phys.}, 18:199--208, 1967.

\bibitem{DaPratoPADOVA1980}
G.~Da~Prato and M.~Iannelli.
\newblock Linear integro-differential equations in {B}anach spaces.
\newblock {\em Rend. Sem. Mat. Univ. Padova}, 62:207--219, 1980.

\bibitem{FujishiroYAMAMOTO2014}
K.~Fujishiro and M.~Yamamoto.
\newblock Approximate controllability for fractional diffusion equations by
  interior control.
\newblock {\em Appl. Anal.}, 93(9):1793--1810, 2014.

\bibitem{GuerreroESAIM2013}
S.~Guerrero and O.~Y. Imanuvilov.
\newblock Remarks on non controllability of the heat equation with memory.
\newblock {\em ESAIM Control Optim. Calc. Var.}, 19(1):288--300, 2013.

\bibitem{HalanaySCL2012}
A.~Halanay and L.~Pandolfi.
\newblock Lack of controllability of the heat equation with memory.
\newblock {\em Systems Control Lett.}, 61(10):999--1002, 2012.

\bibitem{HalanayJMAA2014}
A.~Halanay and L.~Pandolfi.
\newblock Lack of controllability of thermal systems with memory.
\newblock {\em Evol. Equ. Control Theory}, 3(3):485--497, 2014.

\bibitem{HalanayJMAA2015}
A.~Halanay and L.~Pandolfi.
\newblock Approximate controllability and lack of controllability to zero of
  the heat equation with memory.
\newblock {\em J. Math. Anal. Appl.}, 425(1):194--211, 2015.

\bibitem{HalanayDCDS-A2015}
A.~Halanay and L.~Pandolfi.
\newblock Noncontrollability for the {C}olemann-{G}urtin model in several
  dimensions.
\newblock {\em Discrete Contin. Dyn. Syst.}, pages 588--595, 2015.

\bibitem{KATOperturbTheoryBOOK}
T.~Kato.
\newblock {\em Perturbation theory for linear operators}.
\newblock Springer-Verlag, Berlin-New York, second edition, 1976.
\newblock Grundlehren der Mathematischen Wissenschaften, Band 132.

\bibitem{WarmaDISCRCONTDYNSYS2016}
V.~Keyantuo and M.~Warma.
\newblock On the interior approximate controllability for fractional wave
  equations.
\newblock {\em Discrete Contin. Dyn. Syst.}, 36(7):3719--3739, 2016.

\bibitem{LasieckaTriggianiLIBROencVOL1}
I.~Lasiecka and R.~Triggiani.
\newblock {\em Control theory for partial differential equations: continuous
  and approximation theories. {I}}, volume~74 of {\em Encyclopedia of
  Mathematics and its Applications}.
\newblock Cambridge University Press, Cambridge, 2000.
\newblock Abstract parabolic systems.

\bibitem{LavanyaBALACH}
R.~Lavanya and K.~Balachandran.
\newblock Controllability results of linear parabolic integrodifferential
  equations.
\newblock {\em Differential Integral Equations}, 21(9-10):801--819, 2008.

\bibitem{Lebeau95COMMpde}
G.~Lebeau and L.~Robbiano.
\newblock Contr\^{o}le exact de l'\'{e}quation de la chaleur.
\newblock {\em Comm. Partial Differential Equations}, 20(1-2):335--356, 1995.

\bibitem{LIONSlibro68}
J.-L. Lions.
\newblock {\em Contr\^{o}le optimal de syst\`emes gouvern\'{e}s par des
  \'{e}quations aux d\'{e}riv\'{e}es partielles}.
\newblock Avant propos de P. Lelong. Dunod, Paris; Gauthier-Villars, Paris,
  1968.

\bibitem{MAINARDIbook10FRACT}
F.~Mainardi.
\newblock {\em Fractional calculus and waves in linear viscoelasticity}.
\newblock Imperial College Press, London, 2010.
\newblock An introduction to mathematical models.

\bibitem{PANDOLFI1989LAA}
L.~Pandolfi.
\newblock Controllability properties of perturbed distributed parameter
  systems.
\newblock {\em Linear Algebra Appl.}, 122/123/124:525--538, 1989.

\bibitem{PPZ}
L.~Pandolfi, E.~Priola, and J.~Zabczyk.
\newblock Linear operator inequality and null controllability with vanishing
  energy for unbounded control systems.
\newblock {\em SIAM J. Control Optim.}, 51(1):629--659, 2013.

\bibitem{Pazy}
A.~Pazy.
\newblock {\em Semigroups of linear operators and applications to partial
  differential equations}.
\newblock Springer-Verlag, New York, 1983.

\bibitem{PruessLIBRO1993}
J.~Pr\"uss.
\newblock {\em Evolutionary Integral Equations and Applications}.
\newblock Birkh\"auser, Basel, 1993.

\bibitem{TaoGAOcontroZEROjmma16}
Q.~Tao and H.~Gao.
\newblock On the null controllability of heat equation with memory.
\newblock {\em J. Math. Anal. Appl.}, 440(1):1--13, 2016.

\bibitem{TuksnakWeiss}
M.~Tucsnak and G.~Weiss.
\newblock {\em Observation and control for operator semigroups}.
\newblock Birkh\"auser Verlag, Basel, 2009.

\bibitem{WARMAapplANAL2017}
M.~Warma.
\newblock On the approximate controllability from the boundary for fractional
  wave equations.
\newblock {\em Appl. Anal.}, 96(13):2291--2315, 2017.

\end{thebibliography}
  \bibliographystyle{plain}

\end{document}